\input epsf
\documentclass[12pt]{amsart}
\usepackage{latexsym, amsbsy, amsmath, amsfonts, amssymb, amsthm, amscd, epsfig, amscd}
\usepackage[all]{xy}
\newcommand{\Z}{\mathbb Z}
\newcommand{\N}{\mathbb N}
\newcommand{\al}{\alpha}  
\newcommand{\be}{\beta} 
\newcommand{\bkl}{Birman--Ko--Lee} 
\newcommand{\BKL}{$SBKL_n^+ $}
\newcommand{\de}{\delta}
\newcommand{\De}{$\delta$\ }
\newtheorem{Theorem}{Theorem}[section]

\newtheorem{Corollary}{Corollary}[section]
\newtheorem{Proposition}{Proposition}[section]
\newtheorem{Remark}{Remark}[section]
\newcommand{\sg}{\sigma} 
\newcommand{\sip}{\sigma_{i}}

\newcommand{\sjp}{\sigma_{j}} 
\allowdisplaybreaks[2]
\numberwithin{equation}{section}
\numberwithin{figure}{section}
\title{About presentations of braid groups and their generalizations}
\author[Vershinin]{V.~V.~Vershinin}
\address{D\'epartement des Sciences Math\'ematiques,
                                     Universit\'e Montpellier II,
Place Eug\'ene Bataillon,
34095 Montpellier cedex 5, France}
\email{ vershini@math.univ-montp2.fr}
\address{ Sobolev Institute of Mathematics, Novosibirsk, 630090,
Russia }
\email{ versh@math.nsc.ru}
\begin{document}
\begin{abstract}
In the paper we give a survey of rather new notions and results 
which generalize classical ones in the theory of braids.
Among such  notions are various inverse monoids of partial braids.
We also observe presentations different from standard Artin presentation for generalizations of braids. Namely, we consider presentations with small 
number
of generators, Sergiescu graph-presentations and Birman-Ko-Lee presentation.
The work of  V.~V.~Chaynikov on the 
word and conjugacy problems 
for the singular braid monoid in Birman-Ko-Lee generators is described as well.
\end{abstract}
\maketitle
\tableofcontents

\section{Introduction}

The purpose of this paper is to give a survey on some recent notions and results
concerning generalizations  of 
the braids. 

Classical braid groups $Br_n$ can be defined in several ways. Either as a 
set of isotopy classes of 
system of $n$ curves in a three-dimensional space (what is the same as 
 the fundamental group of the configuration space of $n$ points on a plane) or as 
the mapping class group of a disc with $n$ points deleted $D_n$ with its 
boundary fixed, what is equivalent to the subgroup of the braid
automorphisms of the automorphism group of a free group 
 $\operatorname{Aut} F_n $. For the exact definitions we make a reference here to a monograph 
on braid, for example the book of C.~Kassel and V.~Turaev \cite{KT}
or to the previous surveys of the author \cite{Ve3, Ve6, Ve9_5}.
  
The  {\it
pure} braid  group $P_n$ is defined as  the kernel of the 
canonical epimorphism  $\tau_n$ from braids to the symmetric group $\Sigma_n$: 
\begin{equation*}
1\rightarrow P_n\rightarrow Br_n \buildrel
\tau_n \over \longrightarrow \Sigma_n\rightarrow 1 . \label{eq:exact}
\end{equation*}

We fix the canonical Artin presentation \cite{Art1} of the braid group 
 $Br_n$. It has generators $\sigma_i$, 
$i=1, ..., n-1$ and two types of relations: 
\begin{equation}
 \begin{cases} \sigma_i \sigma_j &=\sigma_j \, \sigma_i, \ \
\text{if} \ \ |i-j|
>1,
\\ \sigma_i \sigma_{i+1} \sigma_i &= \sigma_{i+1} \sigma_i \sigma_{i+1}.
\end{cases} \label{eq:brelations}
\end{equation}
The generators $\sigma_i$ correspond to the following automorphisms of $F_n$:
\begin{equation} \begin{cases} 
x_i &\mapsto x_{i+1},
\\ x_{i+1} &\mapsto x_{i+1}^{-1}x_ix_{i+1}, \\
x_j &\mapsto x_j, j\not=i,i+1. 
\end{cases} \label{eq:autf}
\end{equation}

Of course, there exist other presentations of the braid group.
Let 
\begin{equation}
\sigma = \sigma_1 \sigma_{2} \dots \sigma_{n-1},
\label{eq:sigma}
\end{equation} 
then the group $Br_n$ is generated by $\sigma_1$ and $\sigma$ because
\begin{equation}
\sigma_{i+1} =\sigma^i \sigma_1 \sigma^{-i}, \quad i =1, \dots
{n-2}.
\label{eq:sigma_i}
\end{equation} 
The relations for the generators $\sigma_1$ and $\sigma$ are the 
following
\begin{equation}
 \begin{cases}
\sigma_1 \sigma^i \sigma_1 \sigma^{-i} &= 
\sigma^i \sigma_1 \sigma^{-i} \sigma_1 \ \  \text{for} \ \
2 \leq i\leq {n / 2}, \\
\sigma^n &= (\sigma \sigma_1)^{n-1}.
\end{cases} \label{eq:2relations}
\end{equation}
The presentation (\ref{eq:2relations}) was given by Artin in the initial 
paper \cite{Art1}.
This presentation was also mentioned in the books by F.~Klein \cite{Kl}
and by H.~S.~M.~Coxeter and W.~O.~J.~Moser \cite{CM}.

V.~Ya.~Lin in \cite{Li3} gives a slightly different form of this presentation.
Let $\beta\in Br_n$ be defined by the formula
$$\beta= \sigma\sigma_1.$$
Then there is the presentation of the 
group $Br_n$ with generators $\sigma_1$ and $\beta$ and 
 relations:
\begin{equation*}
 \begin{cases}
\beta \sigma^{i-1} \beta &= 
\sigma^i \beta \sigma^{-i-1} \beta \sigma^i \ \  \text{for} \ \
2 \leq i\leq {n / 2}, \\
\sigma^n &= \beta^{n-1}.
\end{cases} \label{eq:2relin}
\end{equation*}
This presentation is called {\it special} in \cite{Li3}.

\begin{figure}[h] 
  \hspace{20pt}\psfig{figure=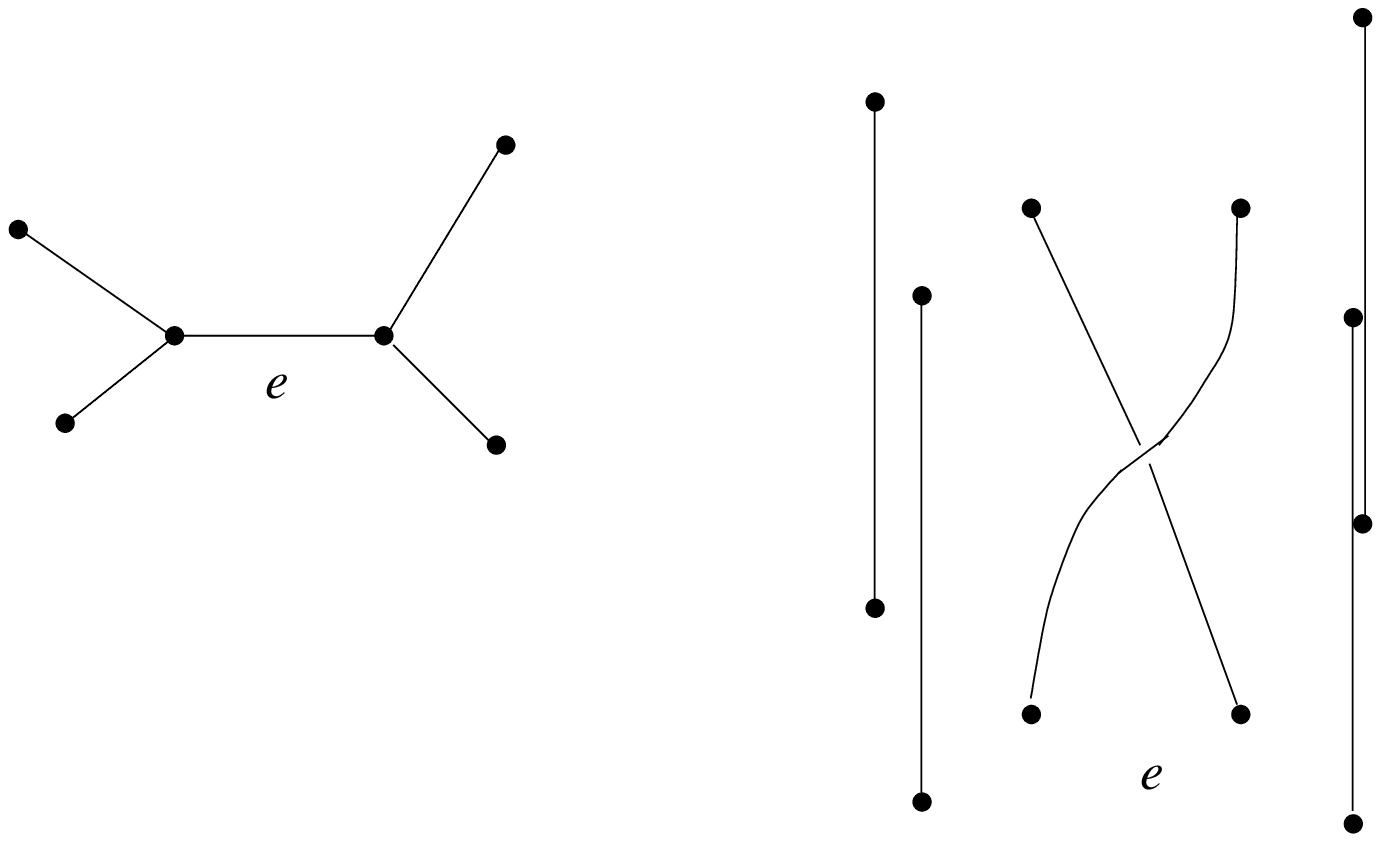,width=8cm}  
  \caption{Edges and geometric braids.} 
  \label{fig:edges}
  \end{figure}

An interesting series of presentations was given by V.~Sergiescu
\cite{Ser}. For every planar graph he constructed a presentation of the
group $Br_n$, where $n$ is the number of vertices of the graph,
with generators corresponding to edges and relations reflecting
the geometry of the graph. 
To each edge $e$ of  the graph he associates the braid $\beta_e$ which is a 
clockwise half-twist along $e$ (see Figure~\ref{fig:edges}). 
Artin's classical presentation (\ref{eq:brelations}) in this context
corresponds to the graph consisting of the interval from 1 to $n$
with the natural numbers (from 1 to $n$) as vertices and with
segments between them as edges. 

To be precise let  $\Gamma$  be a planar  graph. We  call it  \emph{normal}  if  
$\Gamma$  is connected,  and it has no loops or intersections. 
Let $S(\Gamma)$ be  the set of vertices of $\Gamma$. 
If $\Gamma$ is not a tree then we define next what is a 
{\it pseudocycle} on it. The bounded part of the complement
of $\Gamma$ in the plane is the disjoint union of  
a finite number of  open disks $D_1, \dots, D_m$, $m >1$. 
The boundary of  $D_j$ on the plane is a   subgraph $\Gamma(D_j)$ 
of $\Gamma$. 
We choose a point $O$ in the interior of $D_j$,   
and    an edge $\sigma$ of $\Gamma(D_j)$ with vertices  $v_1, v_2$. We suppose 
that the triangle $O v_1 v_2$ is oriented anticlockwise.
We denote $\sg$ by $\sigma(e_{1})$.
We define the {\it pseudocycle associated to} $D_j$ to be the sequence 
of edges $\sg(e_1) \dots \sg(e_{p})$ such that:

-if the vertex $v_{j+1}$ is not  uni-valent, then $\sigma(e_{j+1})$ is the first 
edge on the left of 
$\sigma(e_{j})$ (we consider  $\sigma(e_{j})$ going from  $v_{j}$
to  $v_{j+1}$)  and the vertex $v_{j+2}$ is the other vertex adjacent to 
$\sigma(e_{j+1})$;

-if the vertex $v_{j+1}$ is  uni-valent, then $\sigma(e_{j+1})=\sigma(e_j)$
and  $v_{j+2}=v_{j}$.

-the vertex $v_{p+1}$ is the vertex $v_1$.

Let $\gamma=\sg(e_1) \dots \sg(e_p)$ be a pseudocycle of  $\Gamma$. 
Let  $i=1, \dots, p$. If  $\sg(e_i)=\sg(e_j)$ for some $j \not=i$, then   
we say  that  
\begin{itemize} 
\item $\sg(e_i)$ is the {\it start edge} of a reverse
if $j = i+1$ (we set $e_{p+1}=e_1$),
\item $\sg(e_i)$ is the 
{\it end edge} of a reverse 
if $j = i-1$ (we set  $e_{0}=e_p$). 
\end{itemize} 

In the following we set $\sg_1 \dots \sg_p$ for the pseudocycle 
$\sg(e_1) \dots \sg(e_p)$. 

\begin{figure}[h] 
  \hspace{20pt}\psfig{figure=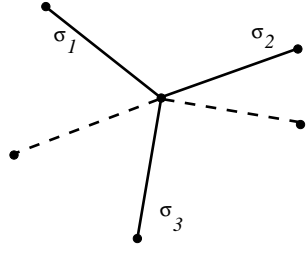,width=4cm}  
  \caption{Nodal relation.} \label{nd:fig}
  \end{figure}

\begin{figure}[h] 
  \hspace{40pt}\psfig{figure=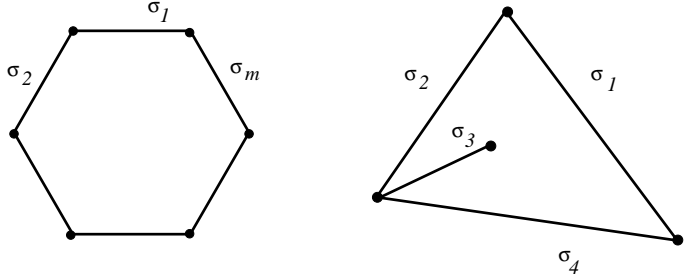,width=9cm}  
  \caption{Pseudocycle relation; on the left $\sg_1 \sg_2 \cdots \sg_{m-1}=
\sg_2 \cdots \sg_{m}=\dots= \sg_{m} \cdots \sg_{m-2}$. 
On the right  $\sg_1 \sg_2 \sg_3^2=  \sg_2 \sg_3^2 \sg_4=\sg_3^2 \sg_4 \sg_1$
and  
$\sg_3 \sg_4 \sg_1\sg_2= \sg_4 \sg_1\sg_2 \sg_3$.}   \label{ps:fig}
\end{figure}

\begin{Theorem} {\rm (V.~Sergiescu \cite{Ser})}  \label{cor:sf1}  
Let  $\Gamma$ be  a normal planar graph   with  $n$ vertices. The braid group
$Br_n$ admits a presentation 
$\langle X_\Gamma \, | \,  R_\Gamma \rangle$, where  $X_\Gamma=\{\sigma \; |\; 
\sigma \; \mbox{is an edge of } \; \Gamma\} $ 
and  $R_\Gamma$ is the set of  following relations: 
\begin{itemize} 
        \item Disjointedness relations (DR): if $\sip$ and $\sjp$ are disjoint, 
        then $\sip \sjp =  \sjp \sip$; 
        \item Adjacency relations (AR): if $\sip, \sjp$ have a common vertex,
then 
        $\sip \sjp \sip $ $= \sjp \sip \sjp$; 
        \item Nodal relations (NR): if  $\{ \sigma_{1}, \sigma_{2}, 
\sigma_{3}\}$  
        have only one common vertex and they are clockwise oriented 
(Figure \ref{nd:fig}), then  
                $$ 
                \sigma_{1} \sigma_{2} \sigma_{3}\sigma_{1}= \sigma_{2}
\sigma_{3} 
                \sigma_{1}\sigma_{2}\,; 
                $$ 
        \item   Pseudocycle relations  (PR): if $\sg_1 \dots \sg_m$ 
     is  a  pseudocycle  and  $\sg_1$ is not the start 
edge or  $\sg_m$ the end edge 
        of a reverse (Figure \ref{ps:fig}), then  
                $$ 
                \sigma_{1} \sigma_{2} \cdots \sigma_{m-1}= \sigma_{2} \sigma_{3} 
\cdots 
                \sigma_{m}\, . 
                $$ 
                   \end{itemize} \label{Theorem:serg}
                  \end{Theorem}

\begin{Remark}  
Theorem \ref{cor:sf1} is true for infinite graphs. Let $\Gamma$ be the direct limit
of its finite subgraphs $\Gamma_i$, then the braid group $Br_\Gamma$ is the  
direct limit of the subgroups $Br_{\Gamma_i}$. 
\end{Remark} 

The graph presentation of Sergiescu underlines the geometric character of 
braids, its connection with configuration spaces. In this survey we confirm this
proposing a thesis: for every generalization of braids of {\it geometric}
character there exists a graph presentation.
  
\begin{figure}
\epsfbox{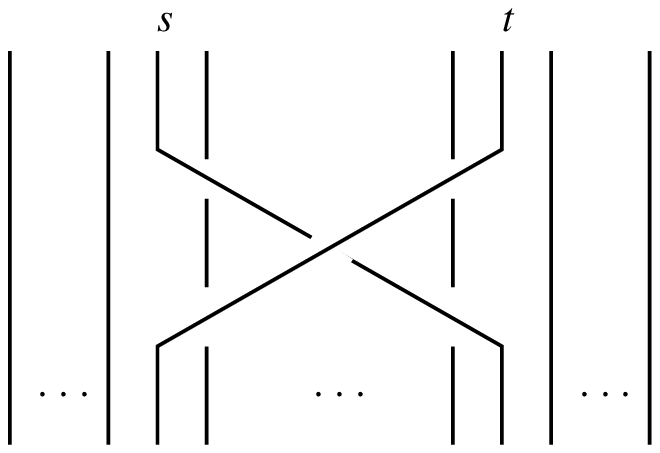}
\caption{}
\label{fi:sbige2}
\end{figure}

Birman,  Ko and Lee \cite{BKL} introduced the presentation 
with the generators $a_{ts}$ with $1 \leq s<t\leq n$
and relations
\begin{equation*} \begin{cases}
a_{ts}a_{rq}&=a_{rq}a_{ts} \ \ {\rm for} \ \ (t-r)(t-q)(s-r)(s-q)>0,\\ 
a_{ts}a_{sr} &=a_{tr}a_{ts}=a_{sr}a_{tr}  \ \ {\rm for} \ \ 
1\leq r<s<t\leq n .
\end{cases}\label{eq:rebkl}
\end{equation*}
The generators $a_{ts}$ are expressed by the canonical generators 
$\sigma_i$
in the following form:
 \begin{equation*} a_{ts}=(\sigma_{t-1}\sigma_{t-2}\cdots\sigma_{s+1})\sigma_s
(\sigma^{-1}_{s+1}\cdots\sigma^{-1}_{t-2}\sigma^{-1}_{t-1})  \ \ 
{\rm for} \ \ 1\leq s<t\leq n.
\label{eq:ats}
\end{equation*} 
Geometrically the generators $a_{s,t}$ are depicted in Figure~\ref{fi:sbige2}.
These generators are very natural and for this presentation 
Birman, Ko and Lee proposed an algorithm which solves the
word problem with the speed $\mathcal O(m^2n)$ while Garside algorithm \cite{Gar}
improved by W.~Turston has a speed $\mathcal O(m^2n\log n)$, where $m$
is the length of a word and $n$ is the number of strands
(see \cite{E_Th}, Corollary~9.5.3). The question of generalization of this 
presentation for other types of braids was raised
in \cite{BKL}.

In Section~\ref{sec:gen}
we describe generalizations of braids that will be involved. In Section~\ref{sec:fewg}
we give the presentations with few generators, in Section~\ref{sec:graphp}
we study graph-presentations in the sense of V.~Sergiescu and in Section~\ref{sec:BKLp}
we give the Birman-Ko-Lee presentation for the singular braid monoid.  
In Section~\ref{sec:chai} we describe the work of 
V.~V.~Chaynikov \cite{Ch} on the 
word and conjugacy problems 
for the singular braid monoid in Birman-Ko-Lee generators.
In Sections~\ref{sec:invmon} -- \ref{sec:m_gen_br} we study inverse 
monoids of partial braids. 

The author is thankful to the organizers of Knots in Poland III
J\'ozef Przytycki and Pawel Trazcyk for the excellent conference.

\section{Generalizations of Braids\label{sec:gen}}

It is interesting to obtain the analogues of the presentations 
mentioned in the Introduction  for various generalizations of braids
 \cite{Bae}, \cite{Bir2}, \cite{Bri1}, \cite{Del}, \cite{FRR2},
\cite{Ve4}.

\subsection{Artin-Brieskorn braid groups\label{subsection:abbg}}

Let  $I$  be a set and  $M =(m_{i,j})$, is a matrix, $m_{i,j}\ \in \N^+\cup \{\infty\}$, $i, j \in I$ with 
the following conditions:  $m_{i,i}= 1$ and
$m_{i,j} >1 $ for $i\not= j$.
J.~Tits in \cite{T} defines the {\it  Coxeter group of  type} $M$
as a group with generators $w_i$, $i\in I$ and relations 
\begin{equation*}
 (w_i w_j)^{m_{i,j}}=e, \ i, j \in I.
 \label{cox_rel}
\end{equation*} 
The corresponding braid groups, which are called {\it Artin-Tits
groups} have the elements $s_i$, $i\in I$ as the 
generators and 
the following set of defining relations:
\begin{equation*}\operatorname{prod}
(m_{i,j};s_i,s_j)= \operatorname{prod} (m_{j,i};s_j,s_i),
\label{ti_br_rel}
\end{equation*} 
where $\operatorname{prod}(m;x,y)$ denotes the product
$xyxy... $ ($m$ factors).

Classification of  irreducible finite
Coxeter groups is well known (see for example Theorem 1, Chapter
VI, \S 4 of \cite{Bo}). It consists of the three infinite series: $A$,
$B$  and $D$ as well as the
exceptional groups $ E_6, E_7, E_8, F_4$, $G_2$, $H_3$, $H_4$ and
$I_2(p).$ 

Let $N$ be a finite set of cardinality $n$, say $N= \{v_1, \dots, v_n\}$.
 Let us equip elements of
$N$ with the signs, i.e. let $SN = \{\delta_1 v_1, \dots, \delta_n v_n\}$,
where $\delta_i =\pm 1$. The Coxeter group $W(B_n)$ of type $B$ can be 
interpreted as 
a group of signed permutations of the set $SN$: 
\begin{equation}
W(B_n)=\{\sigma - \text{
bijection \ of } SN: (-x)\sigma  =-(x)\sigma \text{ for} \ x\in SN\}. 
\label{weylb}
\end{equation}

  The {\it generalized braid group} (or 
{\it Artin--Brieskorn group}) $Br(W)$
{\it of} $W$ \cite{Bri1}, \cite{Del} correspond to the case of 
finite Coxeter group $W$.
 The classical braids on $k$ strings $Br_k$
are obtained by this construction if $W$ is the symmetric group on
$k$ symbols. In this case $m_{i,i+1}=3$, and $m_{i,j}=2$ if $j
\neq i, i+1 $.

The braid group of type $B_n$ has the canonical presentation
with generators  $\sigma_i$, $i=1, \dots, n-1$ and $\tau$, and relations:
\begin{equation}
\begin{cases}
\sigma_i\sigma_j &=\sigma_j\sigma_i, \ \
\text{if} \ \ |i-j| >1,
\\ \sigma_i \sigma_{i+1} \sigma_i &=
\sigma_{i+1} \sigma_i \sigma_{i+1},\\
\tau\sigma_i &=\sigma_i\tau, \ \ \text{if} \ \  i\geq 2,\\
\tau\sigma_1\tau\sigma_1&=\sigma_1\tau\sigma_1\tau . \\
\end{cases} \label{eq:relB}
\end{equation}

\begin{figure}[h] 
\hspace{30pt}\psfig{figure=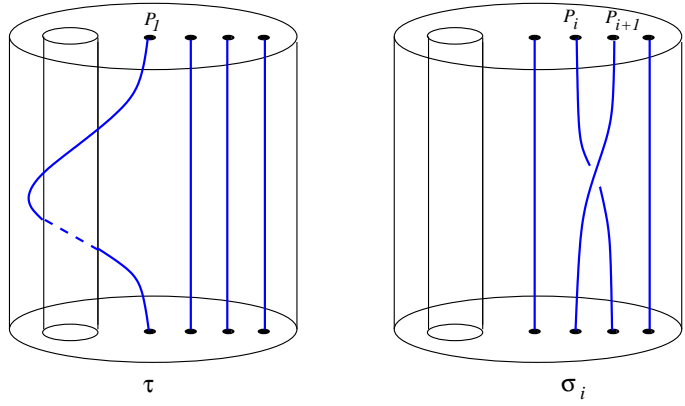,width=9cm}  
\caption{Geometric interpretation of generators $\tau, \sigma_1, \dots, 
\sigma_{n-1}$ of $Br_n(Ann)$.} 
\label{genann:fig}
\end{figure}

This group can be identified with  
the fundamental group of the configuration space of distinct points on 
the plane with one point deleted 
\cite{La}, \cite{Ve0}, what is the same as the braid group on $n$ strands on the 
annulus, $Br_n(Ann)$. 
A geometric interpretation of generators $\tau, \sigma_1, \dots, \sigma_{n-1}$ 
is given in Figure \ref{genann:fig}.

The braid groups of the type $D_n$
has the canonical presentation
with generators  $\sigma_i$ and $\rho$, and relations:
\begin{equation}
\begin{cases}
\sigma_i\sigma_j &=\sigma_j\sigma_i \ \
\text{if} \ \ |i-j| >1,\\ \sigma_i \sigma_{i+1} \sigma_i &=
\sigma_{i+1} \sigma_i \sigma_{i+1},\\
\rho\sigma_i &=\sigma_i\rho \ \ \text{if} \ \  i =1, 3, \dots, n-1,\\
\rho\sigma_2\rho&=\sigma_2\rho\sigma_2 . \\
\end{cases} \label{eq:relD}
\end{equation}

Let V be a complex finite dimensional vector space. A {\it pseudo-reflection} of $GL(V)$ is a
non trivial element $s$ of $GL(V$) which acts trivially on a hyperplane, called the reflecting
hyperplane of $s$.
Suppose that
$W$ is a
finite subgroup of $GL(V)$ generated by  pseudo-reflections;
the corresponding braid groups were studied by M.~Brou\'e, G.~Malle and 
R.~Rouquier \cite{BMR} and also by D.~Bessis and J.~Michel \cite{BM}. 
As in the classical case these groups  can be defined as fundamental
groups of complement in $V$ of the reflecting hyperplanes. The following
classical conjecture generalizes the case of braid groups:

{\it The universal cover of complement in $V$ of the reflecting hyperplane is contractible.} 

(See for example the book by Orlik and Terao \cite{orlikterao}, p. 163 \& p. 259)

This conjecture was proved  by David Bessis \cite{Be}. 
It means that these groups has naturally defined finite dimensional
manifold as $K(\pi,1)$-spaces.

\subsection{Braid groups on surfaces}
Let $\Sigma$ be surface. The $n$th braid group of $\Sigma$ can be defined as 
 the fundamental group of configuration space of $n$ points on $\Sigma$. 
Let $\Sigma$ be a sphere. The corresponding braid group $Br_n(S^2)$  has 
simple geometric interpretation
as a group of isotopy classes of braids lying in a layer between two 
concentric spheres. It has the 
presentation with generators $\delta_i$, $i=1, ..., n-1$, 
which satisfy the braid relations  (\ref{eq:brelations}) and the 
following sphere relation: 
\begin{equation}
\delta_1 \delta_2 \dots \delta_{n-2}\delta_{n-1}^2\delta_{n-2} \dots
\delta_2\delta_1 =1.
 \label{eq:spherelation}
\end{equation}
This presentation was found by O.~Zariski \cite{Za1} in 1936 and then
rediscovered by E.~Fadell and J.~Van Buskirk \cite{FaV} in 1961.

Presentations of braid groups on all closed surfaces were obtained by
G.~P.~Scott \cite{Sc} and others.

 \subsection{Braid-permutation group\label{bp}}

Let $BP_n$ be the subgroup of $\operatorname{Aut} F_n$, generated
by both sets of the automorphisms $\sigma_i$ of (\ref{eq:autf})
and  $\xi_i$ of the following form: 
\begin{equation}
\begin{cases} x_i &\mapsto x_{i+1}, \\ x_{i+1} &\mapsto
x_i,     \\ x_j &\mapsto x_j, j\not=i,i+1, \end{cases}
\label{eq:perm}
\end{equation}
This is the $n$th {\it braid-permutation group} introduced by R.~Fenn,
R.~Rim\'anyi and C.~Rourke   \cite{FRR2} who gave a presentation of
this group: it consists of the set of generators: 
$\{ \xi_i, \sigma_i, \ \ i=1,2,
..., n-1 \}$ 
such that $\sigma_i$ satisfy the braid relations, $\xi_i$ satisfy the symmetric group relations and both of them the satisfy the following mixed relations:
\begin{equation}
 \begin{cases} \sigma_i \xi_j
&=\xi_j \sigma_i, \ \text {if} \  |i-j| >1,
\\ \xi_i \xi_{i+1} \sigma_i &= \sigma_{i+1} \xi_i \xi_{i+1},
\\ \sigma_i \sigma_{i+1} \xi_i &= \xi_{i+1} \sigma_i \sigma_{i+1}.
\end{cases} \label{eq:mixperm}
\end{equation}
\vglue0.01cm
\centerline {The mixed relations for the braid-permutation group}
\smallskip
R.~Fenn, R.~Rim\'anyi and C.~Rourke gave a geometric
interpretation of $BP_n$ as a group of {\it welded braids}.

This group was also studied by A.~G.~Savushkina \cite{Sav2} under the name 
of {\it group 
of conjugating automorphisms} and notation ${\mathbf{C}}_n$. 

Braid-permutation group has an interesting geometric interpretation as a 
motion group. This group was introduced in the PhD thesis of David Dahm, a student of Ralph Fox. It appeared in literature in the paper of Deborah
Goldsmith  \cite{Gol} and then studied by various authors, see, \cite{JMM}, 
for instance.
This is an analogue of the interpretation of the
classical braid group as a mapping class group of a punctured disc. Instead of $n$ points in a disc we consider
$n$ unlinked unknotted circles in a 3-ball. The fundamental group of the
complement of $n$ circles is also the free group $F_n$. Interchanging of two neighbour
points in the case of the braid group corresponds to an automorphism
(\ref{eq:autf}) of the free group. In the case of circles this automorphism
corresponds to a motion of two neighbour circles when one of the circles is
passing inside another one. Simple interchange of two neighbour circles corresponds to the automorphism (\ref{eq:perm}).
 
 Another motivation for studing  braid-permutation groups is given by 
 the {\it  pure braid-permutation group} 
 $ P\varSigma_n$, the kernel of the canonical epimorphism
 $BP_n \to \Sigma_n$. In the context of the motion group it is called
 as the {\it group of loops}, but it has even a longer history and is
 connected with classical works of J.~Nielsen \cite{N} and W.~Magnus
 \cite{magnus} (see also \cite{mks}), as follows.
 Let us denote the kernel of the natural map
$$Aut(F_n) \to GL(n, \mathbb Z)$$ 
by $IA_n$. These groups are similar to the Torelli subgroups of the 
mapping class groups. Nielsen, and
Magnus gave automorphisms which generate $IA_n$ as a group. These 
automorphisms are named as follows:

\begin{itemize}
        \item $\chi_{k,i}$ for $ i \neq k$ with $1 \leq i,k \leq n$, and
        \item $\theta (k;[s,t])$ for $k,s,t$ distinct integers with $1 \leq k,s,t \leq n$ and $s<t$.
\end{itemize}
The definition of the map $\chi_{k,i}$ is given by the formula
\[
\chi_{k,i}(x_j)=
\begin{cases}
x_j & \text{if $k \neq j$,}\\
(x_i^{-1})(x_k)(x_i) & \text{if $k = j$.}
\end{cases}
\] 
The map $\theta(k;[s,t])$ is defined by the formula
\[
\theta(k;[s,t])(x_j)=
\begin{cases}
x_j & \text{if $k \neq j$,}\\
(x_k)\cdot([x_s,x_t]) & \text{if $k = j$.}
\end{cases}
\] for which the commutator is given by 
$[a,b] = a^{-1}\cdot b^{-1}\cdot a \cdot b$. 
 
 The group $IA_2$ is isomorphic to the group of inner
automorphisms ${\rm Inn}(F_2)$, which is isomorphic to the free
group $F_2$. The group $IA_3$ is not finitely presented \cite{km}.

Consider the subgroup of $IA_n$ generated by the $\chi_{k,i}$, the
{\it group of basis conjugating automorphisms} of a free group. This is exactly 
$ P\varSigma_n$. McCool gave a  presentation for it
\cite{mc}. 

The cohomology of $P\varSigma_n$ was computed by C.~Jensen,
J.~McCammond, and J.~Meier \cite{JMM}.  N.~Kawazumi \cite{ka},
T.~Sakasai \cite{s}, T.~Satoh \cite{sa} and A.~Pettet \cite{pe} have
given related cohomological information for $IA_n$. The integral
cohomology of the natural direct limit of the groups $Aut(F_n)$ is
given in work of S.~Galatius \cite{g}.

 \begin{Theorem}{\rm (A.~G.~Savushkina \cite{Sav2})} \label{thm:braid.permutation}  The group $\mathrm{BP}_n$
is the semi-direct product of the symmetric group on $n$-letters
$\Sigma_n$ and the group $P\varSigma_n$ with a split extension

\[
\begin{CD}
1 @>{}>>  P\varSigma_n @>{}>> \mathrm{BP}_n  @>{}>> \Sigma_n @>{}>>
1.
\end{CD}
\]
\end{Theorem}
 The Lie algebra structure obtained from the
descending central series of the group $P\varSigma_n$
was studied by F. R.~Cohen, J. ~Pakianathan, V.~V.~Vershinin and J.~Wu \cite{CPVW} 
and by B.~Berceanu and S.~Papadima \cite{BP}.
Certain subgroups of $P\varSigma_n$ were studied by V.~Bardakov  and R.~Mikhailov \cite{BaMi}.

 \subsection{Singular braid monoid\label{sbm}}
 The set  of singular braids on  $n$ strands,  up to isotopy,  forms  
a  monoid. This is the {\it singular braid monoid} or {\it Baez--Birman monoid} $SB_n$   
\cite{Bae}, \cite{Bir2}. It can be presented as the
monoid with generators $g_i,g_i^{-1},x_i$, $i=1,\dots,n-1,$ and
relations 
\begin{eqnarray*}
&\sigma_i\sigma_j=\sigma_j\sigma_i, \ \text {if} \ \ |i-j| >1,\\
&x_ix_j=x_jx_i, \ \text {if} \ \ |i-j| >1,\\ &x_i\sigma_j=\sigma_jx_i, \
\text {if} \ \ |i-j| \not=1,\\ &\sigma_i \sigma_{i+1} \sigma_i = \sigma_{i+1} \sigma_i
\sigma_{i+1},\\ &\sigma_i \sigma_{i+1} x_i = x_{i+1} \sigma_i \sigma_{i+1},\\ 
&\sigma_{i+1} \sigma_i
x_{i+1} = x_i \sigma_{i+1} \sigma_i,\\ &\sigma_i\sigma_i^{-1}=\sigma_i^{-1}\sigma_i =1.
\end{eqnarray*}
In pictures $\sigma_i$ corresponds to canonical generator of the braid
group and $x_i$ represents an intersection
of the $i$th and $(i+1)$th strand as in 
Figure~\ref{fi:singen}.
The singular braid monoid on two strings is isomorphic to $\Z\oplus\Z^+$. 
 \begin{figure}
\begin{picture}(0,130)(0,-10) 
\put(0,50){\circle*{5}} \put(-100,100){\line(0,-1){100}}
\put(-50,100){\line(0,-1){100}} \put(-25,100){\line(1,-2){50}}
\put(25,100){\line(-1,-2){50}} \put(50,100){\line(0,-1){100}}
\put(100,100){\line(0,-1){100}}
\put(-100,110){\makebox(0,0)[cc]{$1$}}
\put(-50,110){\makebox(0,0)[cc]{$i-1$}}
\put(-25,110){\makebox(0,0)[cc]{$i$}}
\put(25,110){\makebox(0,0)[cc]{$i+1$}}
\put(50,110){\makebox(0,0)[cc]{$i+2$}}
\put(100,110){\makebox(0,0)[cc]{$n$}}
\put(-75,50){\makebox(0,0)[cc]{.\quad.\quad.}}
\put(75,50){\makebox(0,0)[cc]{.\quad.\quad.}}
\end{picture}
\caption{}\label{fi:singen}
\end{figure} 
This monoid embeds in a group $SG_n$ \cite{FKR}
which is called the {\it singular braid group}: $$SB_n\to
SG_n.$$  So, in $SG_n$ the elements $x_i$ become invertible and all
relations of $SB_n$ remain true.

Principal motivations for study of the singular braid monoid lie in the
Vassiliev theory of finite type invariants \cite{Va}. Essential step
in this theory is that a link invariant is extended from usual links to
singular ones. Singular links and singular braids are connected via
singular versions of Alexander theorem proved by Birman \cite{Bir2} and Markov theorem proved by B.~Gemein \cite{Ge1},
so that as well as in the classical case a singular link is an equivalence
class (by conjugation and stabilization) of singular braids. Therefore 
the study of
singular braid monoid especially such questions as conjugation problem
is interesting not only because of its general importance in Algebra but
because of the connections with Knot Theory.

\subsection{Other generalizations of braids that are not considered in the 
paper}

Garside's solution of the word and conjugacy problems for braids had a great 
influence for
the subsequent research on braids. Tools developed by Garside were put as 
the definitions for Gaussian and Garside groups \cite{DP}, \cite{Dehn} or even 
Garside groupoids \cite{Kr}. The later notion is connected also with the 
mapping class groups. 
 
Another direction of generalizations are the parenthesized braids \cite{GrSer},
\cite{Brin}, \cite{Dehn2}. Motivations for these studies are on the other hand 
in D.~Bar-Natan's works on noncommutative tangles \cite{BN1}, \cite{BN2}
and on the other hand in connections with Thompson's group \cite{CFP}.

\section{Presentations of generalizations of braids with few
generators\label{sec:fewg}}

The presentation with two generators  gives an economic way (from the point of view of generators)
to have a vision of the braid group.
We give here the extension of this presentation for the 
natural generalizations of braids.
The results of this section were obtained in \cite{Ve9}.

\subsection{Artin-Brieskorn groups and complex reflexion groups}

For the braid groups of type $B_n$
from the canonical presentation (\ref{eq:relB})
we obtain the presentation with three generators
$\sigma_1$, $\sigma$ and $\tau$ and the following relations:
\begin{equation}
 \begin{cases}
\sigma_1 \sigma^i \sigma_1 \sigma^{-i} &=
\sigma^i \sigma_1 \sigma^{-i} \sigma_1 \ \  \text{for} \ \
2 \leq i\leq {n / 2}, \\
\sigma^n &= (\sigma \sigma_1)^{n-1},\\
\tau\sigma^i\sigma_1\sigma^{-i} &=\sigma^i\sigma_1\sigma^{-i}\tau
\ \ \text{for} \ \ 2 \leq i\leq {n - 2}, \\
\tau\sigma_1\tau\sigma_1&=\sigma_1\tau\sigma_1\tau . \\
\end{cases} \label{eq:2relB}
\end{equation}

If we add the following relations
\begin{equation*}
 \begin{cases}
\sigma_1^2 &= 1, \\
\tau^2 &= 1  \\
\end{cases}
\end{equation*}
to (\ref{eq:2relB}) we then arrive at a presentation of the Coxeter
group of type $B_n$.

Similarly, for the braid groups of the type $D_n$
from the canonical presentation (\ref{eq:relD})
we can obtain the presentation with three generators
$\sigma_1$, $\sigma$ and $\rho$ and the following relations:
\begin{equation}
 \begin{cases}
\sigma_1 \sigma^i \sigma_1 \sigma^{-i} &=
\sigma^i \sigma_1 \sigma^{-i} \sigma_1 \ \  \text{for} \ \
2 \leq i\leq {n / 2}, \\
\sigma^n &= (\sigma \sigma_1)^{n-1},\\
\rho\sigma^i\sigma_1\sigma^{-i} &=\sigma^i\sigma_1\sigma^{-i}\rho
\ \ \text{for} \ \ i=0, 2, \dots,  {n - 2}, \\
\rho\sigma\sigma_1\sigma^{-1}\rho&=\sigma\sigma_1\sigma^{-1}
\rho\sigma\sigma_1\sigma^{-1} . \\
\end{cases} \label{eq:2relD}
\end{equation}

If we add the following relations
\begin{equation*}
 \begin{cases}
\sigma_1^2 &= 1, \\
\rho^2 &= 1  \\
\end{cases}
\end{equation*}
to (\ref{eq:2relD}) we come to a presentation of the Coxeter
group of type $D_n$.

For the exceptional braid groups of types $E_6 - E_8$ our presentations
look similar to the presentation for the groups of type $D$
(\ref{eq:2relD}).
We give it here for $E_8$: it has three generators
$\sigma_1$, $\sigma$ and $\omega$ and the following relations:
\begin{equation}
 \begin{cases}
\sigma_1 \sigma^i \sigma_1 \sigma^{-i} &=
\sigma^i \sigma_1 \sigma^{-i} \sigma_1 \ \  \text{for} \ \
i = 2, 3 , 4, \\
\sigma^8 &= (\sigma \sigma_1)^{7},\\
\omega\sigma^i\sigma_1\sigma^{-i} &=\sigma^i\sigma_1\sigma^{-i}\omega
\ \ \text{for} \ \ i= 0, 1, 3, 4, 5, 6, \\
\omega\sigma^2\sigma_1\sigma^{-2}\omega &=
\sigma^2\sigma_1\sigma^{-2}\omega\sigma^2\sigma_1\sigma^{-2}.
\end{cases} \label{eq:2relE}
\end{equation}
Similarly, if we add the following relations
\begin{equation*}
 \begin{cases}
\sigma_1^2 &= 1, \\
\omega^2 &= 1  \\
\end{cases}
\end{equation*}
to (\ref{eq:2relE}) we arrive at a presentation of the Coxeter
group of type $E_8$.

As for the other exceptional braid groups, $F_4$ has four generators and
it follows from its Coxeter diagram that there is no sense to speak about
analogues of the Artin presentation (\ref{eq:2relations}), $G_2$ and
$I_2(p)$ already have two generators and $H_3$ has three generators. For
$H_4$
it is possible to diminish the number of generators from four to three and
the presentation will be similar to that of $B_4$.

We can summarize informally what we were doing. Let a group have
a presentation
which can be expressed by a ``Coxeter-like" graph. If there exists
a linear subgraph corresponding to the standard presentation of the
classical braid group, then in the ``braid-like" presentation of
our group the part that corresponds to the linear subgraph can be
replaced by two generators and relations (\ref{eq:2relations}).
This recipe can be applied to the complex reflection groups
\cite{ShT} whose  ``Coxeter-like" presentations is obtained in
\cite{BMR}, \cite{BM}. For the series of the complex braid groups
$B(2e, e,r)$, $e\geq 2$, $r\geq 2$ which correspond to the complex
reflection groups $G(de, e, r)$, $d\geq 2$ \cite{BMR} we take the linear
subgraph with nodes $\tau_2, \dots, \tau_r$, and put as above
$\tau = \tau_2 \dots \tau_r$. The group $B(2e, e,r)$ have
presentation with generators $\tau_2, \tau$,
$\sigma, \tau_2^\prime$ and relations
\begin{equation}
 \begin{cases}
\tau_2 \tau^i \tau_2 \tau^{-i} &=
\tau^i \tau_2 \tau^{-i} \tau_2 \ \  \text{for} \ \
2 \leq i\leq {r / 2}, \\
\tau^r &= (\tau \tau_2)^{r-1}, \\
\sigma \tau^i\tau_2\tau^{-i} &= \tau^i\tau_2\tau^{-i}  \sigma,
\ \  \text{for} \ \ 1 \leq i\leq {r - 2}, \\
\sigma \tau_2^\prime\tau_2 &= \tau_2^\prime \tau_2 \sigma, \\
\tau_2^\prime \tau\tau_2\tau^{-1}\tau_2^\prime &=
\tau\tau_2\tau^{-1} \tau_2^\prime\tau\tau_2\tau^{-1}, \\
\tau\tau_2\tau^{-1}\tau_2^\prime \tau_2 \tau\tau_2\tau^{-1}\tau_2^\prime
\tau_2 &= \tau_2^\prime\tau_2 \tau\tau_2\tau^{-1} \tau_2^\prime \tau_2
 \tau\tau_2\tau^{-1}, \\
\underbrace{\tau_2 \sigma \tau_2^\prime\tau_2 \tau_2^\prime
\tau_2 \tau_2^\prime \dots}_{ e+1 \  \text{factors}}&=
 \underbrace{\sigma\tau_2^\prime  \tau_2
\tau_2^\prime \tau_2
\tau_2^\prime \tau_2 \dots}_{ e+1 \  \text{factors}} \ .
\end{cases} \label{eq:2relBde}
\end{equation}
If we add the following relations
\begin{equation*}
 \begin{cases}
\sigma^d &= 1, \\
\tau_2^2 &= 1,  \\
{{\tau_2}^\prime}^2 &= 1  \\
\end{cases}
\end{equation*}
to (\ref{eq:2relBde}) we come to a presentation of the complex
reflection group $G(de, e, r)$.

The braid group $B(d,1, n)$, $d>1$, has the same presentation as the
Artin -- Brieskorn group of type $B_n$, but if we
add the following relations
\begin{equation*}
 \begin{cases}
\sigma_1^2 &= 1, \\
\tau^d &= 1  \\
\end{cases}
\end{equation*}
to (\ref{eq:2relB}) then we arrive at a presentation of the complex
reflection group $G(d, 1, n)$, $d\geq 2$.

For the series of braid groups $B(e, e,r)$, $e\geq 2$, $r\geq 3$
which correspond to the complex
reflection groups $G(e, e, r)$, $e\geq 2$, $r\geq 3$
  we take again the linear
subgraph with the nodes $\tau_2, \dots, \tau_r$, and put as above
$\tau = \tau_2 \dots \tau_r$. The group $B(e, e,r)$ may have
the presentation with generators $\tau_2, \tau$,
$\tau_2^\prime$ and relations
\begin{equation}
 \begin{cases}
\tau_2 \tau^i \tau_2 \tau^{-i} &=
\tau^i \tau_2 \tau^{-i} \tau_2 \ \  \text{for} \ \
2 \leq i\leq {r / 2}, \\
\tau^r &= (\tau \tau_2)^{r-1}, \\
\tau_2^\prime \tau\tau_2\tau^{-1}\tau_2^\prime &=
\tau\tau_2\tau^{-1} \tau_2^\prime\tau\tau_2\tau^{-1}, \\
\tau\tau_2\tau^{-1}\tau_2^\prime \tau_2 \tau\tau_2\tau^{-1}\tau_2^\prime
\tau_2 &= \tau_2^\prime\tau_2 \tau\tau_2\tau^{-1} \tau_2^\prime \tau_2
 \tau\tau_2\tau^{-1}, \\
\underbrace{\tau_2 \tau_2^\prime\tau_2 \tau_2^\prime
\tau_2 \tau_2^\prime \dots}_{ e \  \text{factors}}&=
 \underbrace{\tau_2^\prime  \tau_2
\tau_2^\prime \tau_2
\tau_2^\prime \tau_2 \dots}_{ e \  \text{factors}}.
\end{cases} \label{eq:2relBe}
\end{equation}
If $e=2$ then this is precisely the presentation
for the Artin -- Brieskorn group of type $D_r$ (\ref{eq:2relD}).
If we add the following relations
\begin{equation*}
 \begin{cases}
\tau_2^2 &= 1,  \\
{{\tau_2}^\prime}^2 &= 1  \\
\end{cases}
\end{equation*}
to (\ref{eq:2relBe}), then we obtain a presentation of the complex
reflection group $G(e, e, r)$, $e\geq 2$, $r\geq 3$.

As for the exceptional (complex) braid groups, it is reasonable to
consider the groups $Br(G_{30})$, $Br(G_{33})$ and
 $Br(G_{34})$  which correspond to the complex reflection groups
$G_{30}$, $G_{33}$ and $G_{34}$.

The presentation for $Br(G_{30})$ is similar to the presentation
(\ref{eq:2relB}) of $Br(B_4)$ with the last relation replaced by the
relation of length 5: the three generators
$\sigma_1$, $\sigma$ and $\tau$ and the following relations:
\begin{equation}
 \begin{cases}
\sigma_1 \sigma^2 \sigma_1 \sigma^{-2} &=
\sigma^2 \sigma_1 \sigma^{-2} \sigma_1 , \\
\sigma^4 &= (\sigma \sigma_1)^{3},\\
\tau\sigma^i\sigma_1\sigma^{-i} &=\sigma^i\sigma_1\sigma^{-i}\tau
\ \ \text{for} \ \ i = 2, 3, \\
\tau\sigma_1\tau\sigma_1\tau&=\sigma_1\tau\sigma_1\tau \sigma_1. \\
\end{cases} \label{eq:2relBG30}
\end{equation}
If we add the following relations
\begin{equation*}
 \begin{cases}
\sigma_1^2 &= 1, \\
\tau^2 &= 1  \\
\end{cases}
\end{equation*}
to (\ref{eq:2relBG30}), then we obtain a presentation of complex
reflection group $G_{30}$.

As for the groups $Br(G_{33})$ and $Br(G_{34})$,
we give here the presentation for
the latter one because the ``Coxeter-like" graph for $Br(G_{33})$
has one node less in the linear subgraph (discussed earlier) than that
of $Br(G_{34})$.
This presentation has the three generators
$s$, $z $ ($z=stuvx$ in the reflection generators) and $w$ and the
following relations:
\begin{equation}
 \begin{cases}
s z^i s z^{-i} &=
z^i s z^{-i} s \ \  \text{for} \ \
i = 2, 3 , \\
z^6 &= (z s)^{5},\\
w z^i s w^{-i} &=z^i s z^{-i}w
\ \ \text{for} \ \ i= 0, 3, 4, \\
w z^i s z^{-i} w &=
z^i s z^{-i} w z^i s z^{-i} \ \  \text{for} \ \
i = 1, 2, \\
w z^2 s z^{-2} w z s z^{-1} w z^2 s z^{-2}&=
 z s z^{-1} w z^2 s z^{-2} w z s z^{-1} w . \\
\end{cases} \label{eq:2relBG34}
\end{equation}
In the same way if we add the following relations
\begin{equation*}
 \begin{cases}
s^2 &= 1, \\
w^2 &= 1  \\
\end{cases}
\end{equation*}
to (\ref{eq:2relBG34}), then we come to a presentation of the complex
reflection group $G_{34}$.

We can obtain presentations with few generators for
the other complex reflection groups using the already observed presentations
of the braid groups. For $G_{25}$ and $G_{32}$ we can use the
presentations (\ref{eq:2relations}) for the classical braid groups $Br_4$
and $Br_5$ with the only additional relation
\begin{equation*}
\sigma_1^3 = 1. 
\end{equation*}
\subsection{Sphere braid groups: few generators}  
It has
two generators $\delta_1$, $\delta$ which satisfy relations (\ref{eq:2relations})
(where $\sigma_1$ is replaced by $\delta_1$, and  $\sigma$ is replaced by $\delta$) 
and the following sphere relation:
\begin{equation*}
\delta^n  (\delta_1 \delta^{-1})^{n-1}= 1.
\end{equation*}
\subsection{Braid-permutation groups}
 
For the case of the braid-permutation group $BP_n$ we  add the new
generator
$\sigma$, defined by (\ref{eq:sigma}) to the set of standard generators of
$BP_n$; then relations (\ref{eq:sigma_i}) and the following relations hold
\begin{equation*}
\xi_{i+1} =\sigma^i \xi_1 \sigma^{-i}, \quad i =1, \dots,
{n-2}.
\end{equation*}
This gives a possibility to get rid of $\xi_i$ as well as of $\sigma_i$
for $i\geq 2$.

\begin{Theorem} The braid-permutation group $BP_n$ has a presentation with
generators $\sigma_1$, $\sigma$,  and $\xi_1$ and relations
\begin{equation*}
 \begin{cases}
\sigma_1 \sigma^i \sigma_1 \sigma^{-i} =
\sigma^i \sigma_1 \sigma^{-i} \sigma_1 \ \  \text{for} \ \
2 \leq i\leq {n / 2}, \\
\sigma^n = (\sigma \sigma_1)^{n-1},\\
\xi_1\sigma^i\sigma_1\sigma^{-i}= \sigma^i\sigma_1\sigma^{-i} \xi_1
\ \  \text{for} \ \ i= 2 \dots {n - 2}, \\
\xi_1 \sigma^i \xi_1 \sigma^{-i} =
\sigma^i \xi_1 \sigma^{-i} \xi_1 \ \ \text{for} \ \ i= 2 \dots {n - 2}, \\
\xi_1\sigma\xi_1\sigma^{-1}\sigma_1 = \sigma\sigma_1\sigma^{-1}\xi_1
\sigma \xi_1\sigma^{-1}, \\
\xi_1\sigma\xi_1\sigma^{-1}\xi_1 = \sigma\xi_1\sigma^{-1}\xi_1
\sigma \xi_1\sigma^{-1}, \\
\xi^2=1.
\end{cases}
\end{equation*}     
\end{Theorem}

\subsection{Few generators for the singular braid monoid}
If we add the new generator $\sigma$,
defined by (\ref{eq:sigma}) to the set of generators of $SB_n$ then the following
relations hold
\begin{equation}
x_{i+1} =\sigma^i x_1 \sigma^{-i}, \quad i =1, \dots
{n-2}.
\label{eq:x_i}
\end{equation}
This gives a possibility to get rid of $x_i$, $i\geq 2$.

\begin{Theorem} The singular braid monoid $SB_n$ has a presentation with
generators $\sigma_1$, $\sigma_1^{-1}$, $\sigma$, $\sigma^{-1}$ and
$x_1$
and relations
\begin{equation}
 \begin{cases}
\sigma_1 \sigma^i \sigma_1 \sigma^{-i} =
\sigma^i \sigma_1 \sigma^{-i} \sigma_1 \ \  \text{for} \ \
2 \leq i\leq {n / 2}, \\
\sigma^n = (\sigma \sigma_1)^{n-1},\\
x_1\sigma^i\sigma_1\sigma^{-i}= \sigma^i\sigma_1\sigma^{-i} x_1
\ \  \text{for} \ \ i=0, 2, \dots, {n - 2}, \\
x_1 \sigma^i x_1 \sigma^{-i} =
\sigma^i x_1 \sigma^{-i} x_1 \ \  \text{for} \ \
2 \leq i\leq {n / 2}, \\
\sigma^n x_1 = x_1\sigma^n,\\
x_1\sigma\sigma_1\sigma^{-1}\sigma_1 = \sigma\sigma_1\sigma^{-1}\sigma_1
\sigma x_1\sigma^{-1}, \\
\sigma_1\sigma_1^{-1}=\sigma_1^{-1}\sigma_1 =1,\\
\sigma\sigma^{-1}=\sigma^{-1}\sigma =1.
\end{cases} \label{eq:2gens}
\end{equation}
\end{Theorem}

\section{Graph-presentations \label{sec:graphp}}
\subsection{Braid groups of type $B$ via graphs\label{subsection:bgvg}} 
Graph presentations for the braid groups of the type $B$ and for the singular 
braid monoid 
were studied by the author. 
We recall that the group  $Br_n(Ann)$ embeds in the braid group $Br_{n+1}$ 
as the subgroup of braids with the first strand fixed. 
 
In the following we consider a normal  planar graph $\Gamma$ such that  
there exists a distinguished 
vertex $v$ and such that the graph $\Gamma$ minus the vertex $v$ and all the 
edges adjacent to $v$ is connected also. 
  We call such $\Gamma$ a $1$-{\it punctured} graph. 
 
\begin{figure}[h] 
\hspace{30pt}\psfig{figure=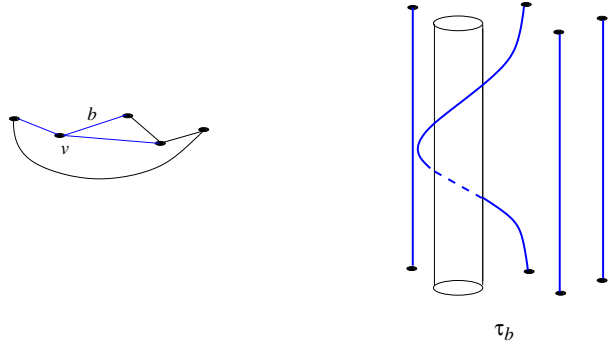,width=8cm}  
\caption{Geometric interpretation of $\tau_b$.} \label{geomint:fig}
\end{figure}

\begin{Theorem} \label{thm:ann} 
Let  $\Gamma$ be  a $1$-punctured  graph   with  $n+1$ vertices. The braid 
group $Br_n(Ann)$ 
admits the presentation 
$\langle X_\Gamma \, | \,  R_\Gamma \rangle$, where  $X_\Gamma=
\{\sigma_a, \tau_b \; |\; 
a$ is an edge of $\Gamma$ not adjacent to the distinguished vertex $v$ and $b$ 
is an edge adjacent to $v \}$ 
and  $R_\Gamma$ is the following set of relations: 
\begin{itemize} 
        \item Disjointedness relations (DR): if the edges $a$ and $c$ 
(respectively $b$ and $c$) are disjoint, 
        then $\sg_a \sg_c =  \sg_c \sg_a$  (respectively $\tau_b \sg_c =  
\sg_c \tau_b$); 
        \item Adjacency relations (AR): if the edges $a$ and $c$ 
(respectively $b$ and $c$) have a common vertex, then 
        $\sg_a \sg_c \sg_a=  \sg_c \sg_a \sg_c$  $(\tau_b \sg_c \tau_b \sg_c=  
\sg_c \tau_b \sg_c \tau_b)$; 
        \item Nodal relations (NR): Let $a,b,c$ be three edges such that they 
                                                have only one common vertex and 
they are clockwise ordered.  
                   If  the edges $a,b,c$ are not adjacent to $v$,  then  
               $$ 
               \sigma_{a} \sigma_{b} \sigma_{c}\sigma_{a}= \sigma_{b} \sigma_{c} 
               \sigma_{a}\sigma_{b}\, ; 
               $$ 
if the edges $a,c$ are not adjacent to $v$ and $b$ is adjacent to $v$, then 
                                                 $$ 
               \sigma_{a} \sigma_{b} \tau_{c}\sigma_{a}= \sigma_{b} \tau_{c} 
               \sigma_{a}\sigma_{b}\, ; 
               $$ 
                                                 $$ 
               \tau_b \sigma_{c} \sigma_{a} \tau_b \sigma_{c}= 
\sigma_{a} \tau_b \sigma_{c} 
               \sigma_{a} \tau_b \, ; 
               $$ 
 
\item   pseudocycle relations  (PR): if the edges $a_1, \dots, a_m$ form 
a pseudocycle, $a_1$ is not the start edge or  $a_m$ the end edge 
of a reverse  and all $a_i$ are not adjacent to $v$, then  
                $$ 
                \sigma_{a_1} \sigma_{a_2} \cdots \sigma_{a_{m-1}}= 
\sigma_{a_2} \sigma_{a_3}  \cdots 
                \sigma_{a_m}\, . 
                $$ 
If $a_1, a_m$ are adjacent to $v$, then 
$$ 
                \tau_{a_1} \sigma_{a_2} \cdots \sigma_{a_{m-1}}= 
\sigma_{a_2} \sigma_{a_3}  \cdots 
                \tau_{a_m}\, . 
                $$ 
 
\end{itemize} 
\end{Theorem} 
\begin{Remark} 
As in Theorem \ref{cor:sf1}, the  nodal relation (NR) implies also the equality   
                                                 $$ 
               \sigma_{a} \sigma_{b} \sigma_{c}\sigma_{a}=  
               \sigma_{b} \sigma_{c} \sigma_{a}\sigma_{b}= \sigma_{c} \sigma_{a} 
               \sigma_{b}\sigma_{c}\, . 
               $$ 
\end{Remark} 
The geometric interpretation of generators is the following. The 
distinguished vertex 
corresponds to the deleted  point of the plane. To any edge $a$ that is not 
adjacent to $v$ 
 we associate the corresponding positive half twist. 
To any edge 
$b$ adjacent to $v$ we associate 
the braid $\tau_b$ as in Figure~\ref{geomint:fig}. 

\begin{Remark}  
This Theorem as well as Theorem \ref{cor:sf1} is true for infinite graphs
via the direct limit arguments.
\end{Remark}  
 
\begin{figure}[h] 
\hspace{20pt}\psfig{figure=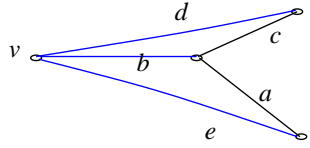,width=4cm}  
\caption{Nodal relation $\tau_b \sigma_{c} \sigma_{a} \tau_b \sigma_{c}= 
\sigma_{a} \tau_b \sigma_{c}\sg_a \tau_b$ 
holds in $B_\Gamma(Ann)$.} \label{nodann:fig}
\end{figure}

\begin{figure}[h] 
\hspace{20pt}\psfig{figure=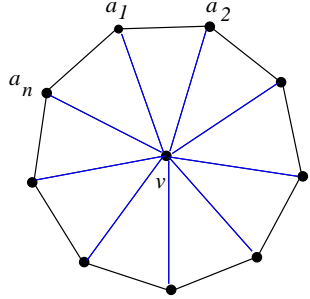,width=4cm}  
\caption{A graph associated to $Br_n(Ann)$.} 
\label{fig:graph}
\end{figure}

\noindent To prove the relation 
$\tau_b \sigma_{c} \sigma_{a} \tau_b \sigma_{c}= 
\sigma_{a} \tau_b \sigma_{c} \sg_a \tau_b$
we add two edges $d$ and $e$, with their corresponding braids 
$\tau_d$ and $\tau_e$ as in Figure \ref{nodann:fig}. The braid $\tau_d$ is 
equivalent 
to the braid $\sg_c^{-1} \tau_b \sg_c$ and the braid $\tau_e$ is equivalent 
to the braid $\sg_a \tau_b \sg_a^{-1}$. Then the braids $\sg_c^{-1} \tau_b
\sg_c$  
and $\sg_a$ commute, as well as $\sg_a \tau_b \sg_a^{-1}$ and $\sg_c$. 
So we have the following equalities, that can be easily verified on 
corresponding braids: 
$$ 
\tau_b \sigma_{c} \sigma_{a} \tau_b \sigma_{c}=  
\sg_c \sg_c^{-1}  \tau_b \sigma_{c} \sigma_{a} \tau_b \sigma_{c}= 
 \sg_c \sg_a \sg_c^{-1}  \tau_b \sigma_{c} \tau_b \sigma_{c}= 
 $$ 
$$ 
 \sg_c \sg_a \sg_c^{-1}  \sg_c \tau_b \sigma_{c} \tau_b= 
  \sg_c {\sg_a \tau_b \sg_a^{-1}} \sg_a \sigma_{c} \tau_b= 
\sg_a \tau_b \sg_a^{-1} \sg_c \sg_a \sigma_{c} \tau_b= 
\sigma_{a} \tau_b \sigma_{c} \sg_a \tau_b \, . 
$$ 
\begin{Corollary} 
The automorphism group of $Br_n(Ann)$ contains a group isomorphic to the  
dihedral group $D_{n-1}$. 
\end{Corollary} 
One can associate to the graph given in Figure \ref{fig:graph} a presentation
for $Br_n(Ann)$. 

It is possible  generalize  such approach to braid groups on a 
planar surface, 
 i.e. a surface of genus $0$ with $l>1$ boundary components.  
 In this case one considers a normal planar graph with $k \,(=l-1)$ distinguished
vertices  
$v_1, \dots, v_k$ such that there are no edges connecting distinguished
vertices and such that the graph $\Gamma$ minus the vertices $v_1, \dots, v_k$  
  and all the edges adjacent to $v_1, \dots, v_k$  is also connected. 
 We label by $\{ \tau_{1,j}, \dots, \tau_{m,j}\}$ the edges adjacent to $v_j$
and by $\{ \sigma_1, \dots, \sigma_p\}$ 
 the edges disjoint from the set $\{ v_1, \dots, v_k \}$.  We say that  $\Gamma$ 
is a $k$-punctured graph. 
 As in Theorem \ref{thm:ann} one can associate to any  $k$-punctured graph 
$\Gamma$ on $n$ vertices 
 a set of generators for the braid group on $n$ strands on surface 
of genus $0$ with $k+1$ boundary components, 
 with the above geometrical interpretation of generators.  

\subsection{Graph-presentations for the surface braid groups}   
These presentations were considered in \cite{BelV}.
Let  $\Gamma$  be a normal graph on an orientable surface $\Sigma$
and $S(\Gamma)$ denotes  the set of vertices of  $\Gamma$. 
In the same way as earlier we associate to the edges of $\Gamma$ the corresponding  
geometric braids on $\Sigma$ (Figure~\ref{fig:edges}) 
and we define  $Br_\Gamma(\Sigma)$ as the subgroup of  $Br_{|S(\Gamma)|}(\Sigma)$ 
generated by these braids. 
\begin{Proposition}\label{prop:ngraph}
Let $\Sigma$ be an oriented surface such that $\pi_1(\Sigma)\not= 1$ and
let $\Gamma$ be a normal graph on $\Sigma$. Then $Br_\Gamma(\Sigma)$ 
is a proper subgroup of $Br_{|S(\Gamma)|}(\Sigma)$.
\end{Proposition}
 
\subsection{Sphere braid groups presentations via graphs}  
 
Now let the surface $\Sigma$ be a sphere $S^2$ and  $\Gamma$ denotes   
a normal finite graph on this sphere.  
We define a pseudocycle as in Introduction:
we consider the set $S^2 \setminus \Gamma$ as the disjoint union of  
a finite number of  open disks $D_1, \dots, D_m$, $m >1$
and define  the {\it pseudocycle associated to} $D_j$ exactly 
in the same way. 

 Let $\Delta$ be a maximal tree of a normal graph $\Gamma$ on $q+1$ vertices.
Then 
$\Delta$ has $q$ edges. Let  $v_1, v_2$ be two vertices adjacent 
to the same edge  $\sg$ of $\Delta$. Write $\sigma(f_{1})$ for $\sg$.
 We define the {\it circuit} $\sg(f_1) \dots \sg(f_{2q})$ as follows:

-if the vertex $v_{j+1}$ is not  uni-valent, then $\sigma(f_{j+1})$ is the 
first edge on the left  of
 $\sigma(f_{j})$ (we consider  $\sigma(f_{j})$ going from  $v_{j}$
to  $v_{j+1}$) and the vertex $v_{j+2}$ is the other vertex adjacent to 
$\sigma(f_{j+1})$;

-if the vertex $v_{j+1}$ is uni-valent, then $\sigma(f_{j+1})=\sigma(f_j)$
and  $v_{j+2}=v_{j}$.

\begin{figure}[h] 
  \hspace{20pt}\psfig{figure=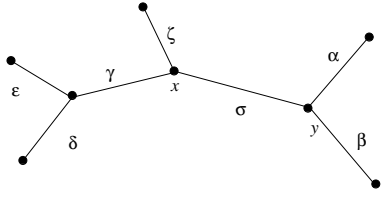,width=5cm}  
  \caption{ $\delta_{x,y} (\Delta) = \sg \al^2 \be^2 \sg \gamma \delta^2
\epsilon^2 \gamma \zeta^2$
and $\delta_{y,x} (\Delta) = \sg \gamma \delta^2 \epsilon^2 \gamma \zeta^2 \sg
\al^2 \be^2  $. } 
\label{fig:deftree}
  \end{figure}

\noindent
This way we come back  to $v_1$ after passing twice through each edge of 
$\Delta$.  Write
$\delta_{v_1,v_2}(\Delta)$ for the word in $X_{\Gamma}$  
corresponding to the  circuit $\sg(f_1) \dots \sg(f_{2q})$ (Figure~\ref{fig:deftree}). 
 
\begin{Theorem}  \label{cor:sf2}  
Let  $\Gamma$ be  a normal  graph   with  $n$ vertices. The braid group
$Br_n(S^2)$ admits a presentation 
$\langle X_\Gamma \, | \,  R_\Gamma \rangle$, where  $X_\Gamma=\{\sigma \; |\; 
\sigma \; \mbox{is an edge of } \; \Gamma\} $ 
and  $R_\Gamma$ is the set of  following relations: 
\par Disjointedness relations (DR);
Nodal relations (NR) (Figure \ref{nd:fig});
Pseudocycle relations  (PR) (Figure \ref{ps:fig}),
                exactly as in Theorem~\ref{Theorem:serg} and the new
\par Tree relations (TR): $\delta_{x, \, y}(\Delta)=1$, 
 for every maximal tree $\Delta$ of $\Gamma$ and  
        every  ordered pair of vertices   $x,y$  such that they are adjacent to the 
same edge  $\sigma$ of $\Delta$. 
\end{Theorem} 

\begin{Remark}  
The statement  of Theorem \ref{cor:sf2} is   highly redundant. 
For instance one can show that a relation (TR)  on a given  
maximal tree of $\Gamma$, together with the relations (DR), (AR), (NR) and (PR),
generate 
the (TR) relation for  any other maximal tree of $\Gamma$. 
Anyway, these presentations are  symmetric and   
one can read off the relations from the geometry of $\Gamma$. 
\end{Remark} 
  
\begin{Remark} \label{rem:thvlad} 
Let $\gamma \subseteq \Gamma$ be a star (a graph which consists of several edges
joined in one point). 
                For any clockwise ordered subset  $\{ \sigma_{i_1}, \dots,$ 
                $  \sigma_{i_j} \, | \, j\ge2 \, \}$ of edges of  $\gamma$ the
following relation holds 
in the group $\langle X_\Gamma \, | \,  R_\Gamma \rangle$: 
                $$ 
                \sigma_{i_1} \dots \sigma_{i_j} \sigma_{i_1}= \sigma_{i_j}
\sigma_{i_1} \dots 
                \sigma_{i_j} \, . 
                $$ 
\end{Remark} 
  
\subsubsection{Geometric interpretation of relations} \label{subsec:geom} 
 
\begin{figure}[h] 
  \hspace{20pt}\psfig{figure=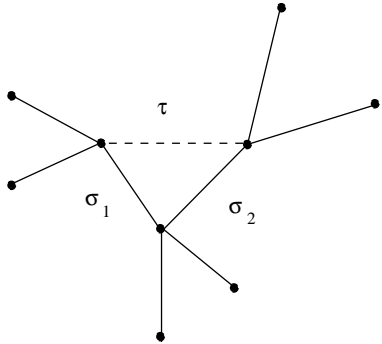,width=5cm}  
  \caption{Adding or removing a triangle.}  \label{fig:addremove}
  \end{figure}

\begin{figure}[h] 
  \hspace{20pt}\psfig{figure=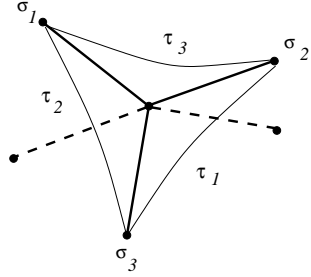,width=4cm}  
  \caption{Nodal relation holds in $B_\Gamma(S^2)$.} \label{fig:nodal}
  \end{figure}

It is geometrically evident that the relations (AR) and (DR) hold  
in $Br_\Gamma(S^2)$. 
Let  $\Gamma$  contain  a triangle $\sg_1,\sg_2, \tau$  as in Figure
\ref{fig:addremove}.  
  Corresponding   braids satisfy the relation 
$\tau = \sg_1 \sg_2 \sg_1^{-1}$ and thus 
$\tau  \sg_1 =\sg_1 \sg_2$ in $Br_\Gamma(S^2)$. 
 The relation $\sg_1 \sg_2 =\sg_2 \tau$ follows from the braid relation 
$\sg_1 \sg_2 \sg_1^{-1}= \sg_2^{-1}\sg_1 \sg_2$. 
Let $\sg_1, \sg_2, \sg_3$ be arranged as in Figure~\ref{fig:nodal}. We 
add three edges $\tau_1, \tau_2, \tau_3$.  
The nodal relation follows  from the pseudocycle relations on triangles 
$\tau_1 \sg_2 \sg_3$, $\tau_2 \sg_1 \sg_3$ and $\tau_3 \sg_1 \sg_2$. 
In fact,  $\sg_1 \sg_2 \sg_3 \sg_1 = \sg_2   \tau_3  \sg_3 \sg_1= 
\sg_2  \sg_3 \tau_3 \sg_1 =  \sg_2  \sg_3 \sg_1 \sg_2$. 
 All other pseudocycle  relations follow from induction on the length of the
cycle. 

\begin{figure}[h] 
  \hspace{20pt}\psfig{figure=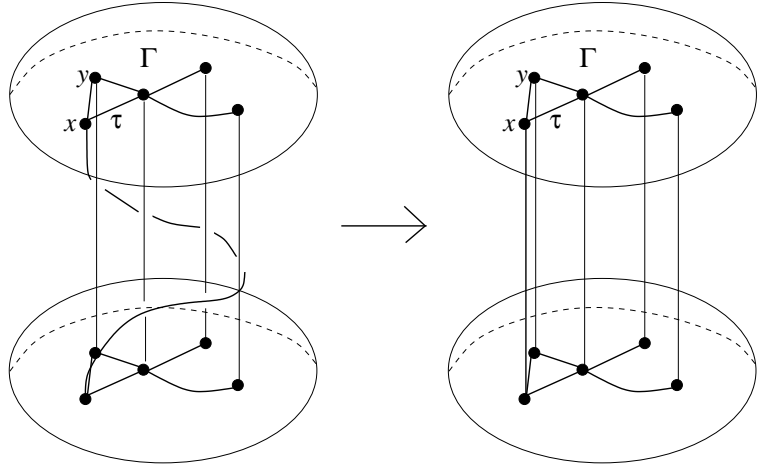,width=10cm}  
  \caption{The braid  $\delta_{x, \sigma} (\Delta)$
 associated to the tree $\Delta= \Gamma \setminus \tau$.}  \label{fig:tretrivial}
  \end{figure}

 Let $\Delta$ be a maximal tree of $\Gamma$. Let  $\sg$ be an edge of $\Delta$ 
and let $x, y$ be the two  adjacent vertices.
The element $\delta_{x, y}(\Delta)$ corresponds to a (pure) braid such that 
the braid obtained by removing  the string starting from the vertex  $x$
is isotopic to the trivial braid. This string  
goes around (with clockwise orientation) all other vertices 
(Figure~\ref{fig:tretrivial} on the left). 
The braid $\delta_{x, y}(\Delta)$ is isotopic to the trivial braid 
in   $Br_\Gamma(S^2)$ and so  $\delta_{x, y}(\Delta) =1$ 
(Figure~\ref{fig:tretrivial}). 
    Therefore the  natural map  $\phi_\Gamma: 
\langle X_\Gamma \, | \,  R_\Gamma\rangle \to Br_\Gamma(S^2)$ is a homomorphism. 

 \subsection{Singular braids and graphs\label{subsection:sinbg}} 
 
 As in the case of classical braids, one can extend the group 
 $Br_n(\Sigma)$ 
 to the monoid $SB_n(\Sigma)$ of singular  braids on $n$ strands on the surface 
$\Sigma$.  
 Presentations for this monoid 
 are given in  \cite{beltes} and \cite{gon3}. 
 
 In this section we provide presentations by graphs for the monoid $SB_n$ 
and for the monoid  
 $SB_n(Ann)$ of singular braids on $n$ strands of the annulus. 
 
 Let $\Gamma$ be a normal planar graph. We associate to any edge $a$ three 
singular braids: 
 $\sg_a$ will denote the positive half-twist associated to $a$ (as in Figure~\ref{fig:edges}), 
$\sg_a^{-1}$  
 will denote  the corresponding negative half-twist and $x_a$ the corresponding 
singular crossing. 
 
 \begin{Theorem} \label{theo:s} 
  Let $\Gamma$ be a normal planar graph with $n$ vertices. 
 The singular braid monoid $SB_n$ has the presentation 
 $\langle X_\Gamma, R_\Gamma \rangle$ where 
 $X_\Gamma = $ 
 $\{\sigma_a, \sigma_a^{-1},$ $ x_a | \ a \ \text{is an edge}$ 
$\text{of} \, \ \Gamma \}$ 
 and $R_\Gamma$ is formed by the following six types of relations: 
 \begin{itemize} 
 \item disjointedness: if the edges  $a$ and $b$ are disjoint, then 
 $$\sigma_a\sigma_b = \sigma_b\sigma_a, \; x_a x_b = x_b x_a, \; 
\sigma_a x_b = x_b\sigma_a,$$ 
 \item commutativity: 
 $$\sigma_a x_a = x_a\sigma_a ,$$ 
 \item invertibility: 
 $$\sigma_a\sigma_a^{-1} = \sigma_a^{-1}\sigma_a =1 ,$$ 
 \item adjacency: if the edges $a$ and $b$ have a common vertex, then 
$$\sigma_a\sigma_b \sigma_a= \sigma_b\sigma_a \sigma_b ,$$ 
$$x_a\sigma_b \sigma_a= \sigma_b\sigma_a x_b ,$$ 
\item nodal: if the edges  $a$, $b$ and $c$  have a common vertex 
and are placed clockwise, then 
$$\sigma_a\sigma_b \sigma_c \sigma_a= \sigma_b\sigma_c \sigma_a 
\sigma_b = \sigma_c\sigma_a \sigma_b \sigma_c , $$ 
$$x_a\sigma_b \sigma_c \sigma_a = \sigma_a\sigma_b \sigma_c x_a , $$ 
$$\sigma_a \sigma_b x_c \sigma_a= \sigma_b x_c \sigma_a \sigma_b , $$ 
$$x_a\sigma_b x_c \sigma_a = \sigma_b x_c \sigma_a x_b , $$ 
\item pseudocycle: if the edges $a_1$, $\dots$, $a_n$ form an 
irreducible pseudocycle and if $a_1$ is not the starting 
edge nor $a_n$ is the end edge of a reverse, then 
$$\sigma_{a_1} \dots \sigma_{a_{n-1}} = \sigma_{a_2}\dots \sigma_{a_n} ,$$ 
$$x_{a_1}\sigma_{a_2} \dots \sigma_{a_{n-1}} = 
\sigma_{a_2}\dots\sigma_{a_{n-1}} x_{a_n} .$$ 
\end{itemize} 
\end{Theorem} 
 
The last aim of this section is to give  graph presentations for 
the singular braid monoid on $n$ strands of the annulus. 
\begin{Theorem}  \label{thm:anns}
The singular braid monoid on $n$ strands of the annulus $SB_n(Ann)$ admits 
the following presentation: 
 
-Generators: 
$\sigma_i, \sigma_i^{-1}, x_i$, $(i=1,\dots,n-1), \tau, \tau^{-1}$. 
\par 
-Relations: 
\begin{eqnarray*} 
&(R1)  &\; \sigma_i\sigma_j=\sigma_j\sigma_i, \; \mbox{if} \quad |i-j| >1 \, ;\\ 
&(R2) & \; x_ix_j=x_jx_i, \; \mbox{if} \quad |i-j| >1 \, ;\\ 
&(R3) & \;x_i\sigma_j=\sigma_j x_i,   \; \mbox{if} \quad |i-j| \not=1 \, ;\\ 
&(R4) & \; \sigma_i \sigma_{i+1} \sigma_i = \sigma_{i+1} \sigma_i \sigma_{i+1}
\, ;\\ 
&(R5) & \; \sigma_i \sigma_{i+1} x_i = x_{i+1} \sigma_i \sigma_{i+1} \, ;\\ 
&(R6) & \; \sigma_{i+1} \sigma_ix_{i+1} = x_i \sigma_{i+1} \sigma_i \, ;\\ 
&(R7) & \; \tau\sigma_1\tau\sigma_1 =  \sigma_1\tau\sigma_1\tau \, ; \\ 
&(R8) & \; \tau\sigma_1\tau x_1 = x_1\tau\sigma_1\tau, \\ 
&(R9) & \; \tau \sigma_i = \sigma_i\tau,  \; \mbox{if} \quad  i \geq 2 \, ;\\ 
&(R10) & \; \tau x_i = x_i\tau,  \; \mbox{if} \quad i \geq 2 \, ;\\ 
&(R11) & \; \sigma_i\sigma_i^{-1} =\sigma_i^{-1}\sigma_i =\tau\tau^{-1} = 
\tau^{-1}\tau =1. 
\end{eqnarray*} 
\end{Theorem} 

The geometric interpretation of $\sg_i$ and $\tau$ is given in Figure 
\ref{genann:fig}. 
 
\begin{figure} 
\hspace{30pt}\psfig{figure=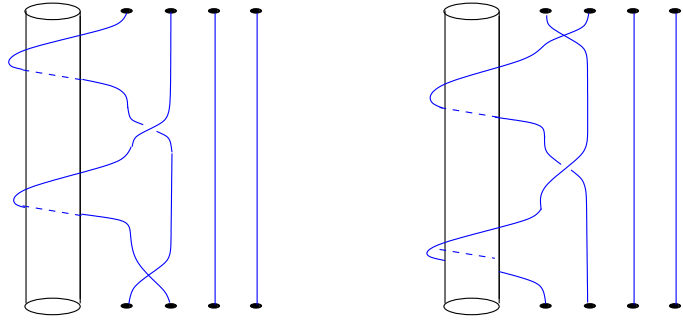,width=9cm}
\caption{The words $\tau\sigma_1\tau x_1$ and $x_1\tau\sigma_1\tau$   
represent the same element in $SB_n(Ann)$.}\label{fi:reltaux} 
\end{figure}

 We get the Reidemeister moves for singular knot theory in 
a solid torus if we add the move depicted on 
Figure~\ref{fi:reltaux} to the regular (without singularities) 
Reidemeister moves of knot theory in a solid torus. 
This Reidemeister move means how a singular point goes around the 
axis of the torus (fixed string). The proof that the list 
R1-R11 is a complete set of relations is standard: every 
isotopy can be decomposed in a sequence of elementary isotopies 
which correspond to relations R1-R11 (see also \cite{gon3}). 
 
\begin{Remark}  The singular braid monoid on $n$ strands of the annulus 
 differs from the singular Artin monoid of type $B$ 
 as defined by R.~Corran \cite{Cor}, where the numbers of singular and 
regular generators are the same. The singular generator associated  to $\tau$
can not be interpreted geometrically as above. 
\end{Remark} 
 
As in Subsection~\ref{subsection:bgvg} we consider 1-punctured graphs. 
To any edge $a$ disjoint from the distinguished vertex $v$ of $\Gamma$ we  
associate three singular braids: 
$\sg_a$ will denote the positive half-twist associated to $a$, 
$\sg_a^{-1}$  
will denote  the corresponding negative half-twist and $\tau_a$ denotes
the corresponding singular crossing. 
 
 The graph presentations for the singular braid monoid in the solid torus 
 arise from Theorems~\ref{theo:s} and~\ref{thm:anns}. 
 \begin{Theorem} Let $\Gamma$ be a one-punctured graph on $n$ vertices. 
The monoid $SB_n(Ann)$ admits the 
presentation $\langle X_\Gamma, R_\Gamma \rangle$, where 

-$X_\Gamma =  \{\sigma_a, \sigma_a^{-1}, x_a, \tau_b, \tau_b^{-1}\}$,
for any  edge $a$ of
$ \, \Gamma $
not incident with the distinguished vertex  $v$, and for any  edge $b$ of 
\, $\Gamma$ 
adjacent to the distinguished vertex 
$v$;

-  $R_\Gamma$ is formed by the relations given in 
Theorems~\ref{thm:ann} and \ref{theo:s} and the following new nodal 
and invertibility relations:
$$\sigma_a\tau_b\sigma_c x_a = x_c\sigma_a\tau_b\sigma_c,$$ 
$$\tau_b\sigma_c\sigma_a\tau_b x_c = x_a\tau_b\sigma_c\sigma_a\tau_b,$$ 
$$\tau_b\tau_b^{-1} = \tau_b^{-1}\tau_b =1.$$ 
\label{theo:x} 
\end{Theorem} 
\section{Birman -- Ko -- Lee presentation for the singular
braid monoid \label{sec:BKLp}}

\begin{figure}
\epsfbox{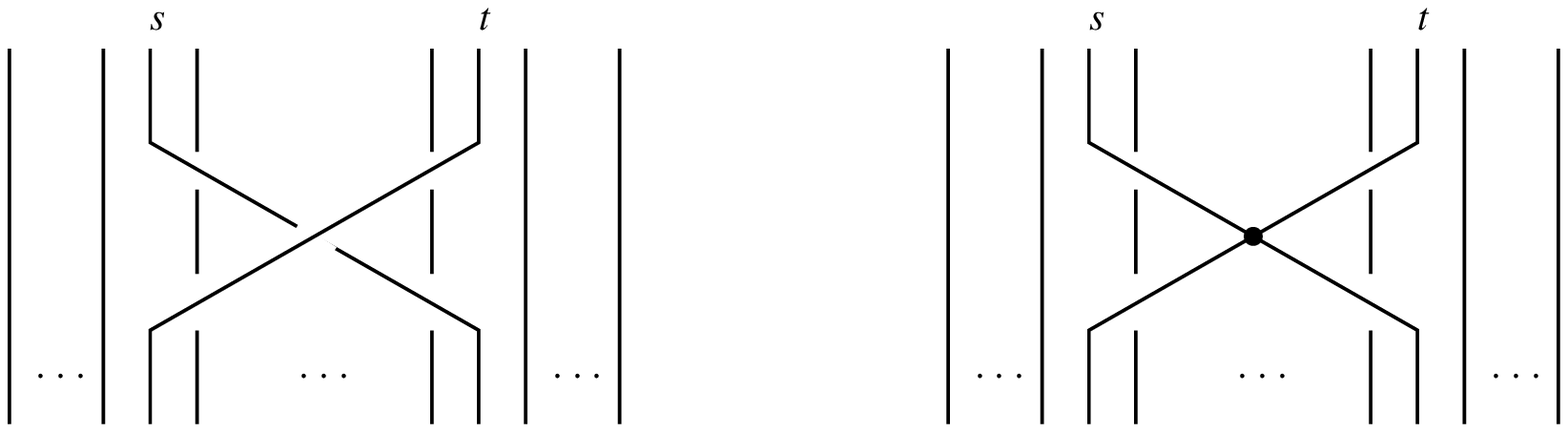}
\caption{}
\label{fi:sbige}
\end{figure}

The analogue of the presentation of Birman,  Ko and Lee
for the singular braid monoid was given in \cite{Ve10}.
For $1\leq s<t\leq n$ and
$1\leq p<q\leq n$ we consider the elements of $SB_n$ which are
defined by
\begin{equation*} \begin{cases}
a_{ts}&=(\sigma_{t-1}\sigma_{t-2}\cdots\sigma_{s+1})\sigma_s
(\sigma^{-1}_{s+1}\cdots\sigma^{-1}_{t-2}\sigma^{-1}_{t-1})  \ \
{\rm for} \ \ 1\leq s<t\leq n, \\
a_{ts}^{-1}&=(\sigma_{t-1}\sigma_{t-2}\cdots\sigma_{s+1})\sigma_s^{-1}
(\sigma^{-1}_{s+1}\cdots\sigma^{-1}_{t-2}\sigma^{-1}_{t-1})  \ \
{\rm for} \ \ 1\leq s<t\leq n, \\
b_{qp}&=(\sigma_{q-1}\sigma_{q-2}\cdots\sigma_{p+1}) x_p
(\sigma^{-1}_{p+1}\cdots\sigma^{-1}_{q-2}\sigma^{-1}_{q-1})  \ \
{\rm for} \ \ 1\leq p<q\leq n.
\end{cases}
\end{equation*}
Geometrically the generators $a_{s,t}$ and $b_{s,t}$ are depicted
in Figure~\ref{fi:sbige}.
\begin{Theorem} The singular braid monoid $SB_n$ has a presentation with
generators $a_{ts}$, $a_{ts}^{-1}$
for $1\leq s<t\leq n$ and  $b_{qp}$ for
$1\leq p<q\leq n$
and relations
\begin{equation} \begin{cases}
a_{ts}a_{rq}&=a_{rq}a_{ts} \ \ {\rm for} \ \ (t-r)(t-q)(s-r)(s-q)>0,\\
a_{ts}a_{sr} &=a_{tr}a_{ts}=a_{sr}a_{tr}  \ \ {\rm for} \ \
1\leq r<s<t\leq n , \\
a_{ts}a_{ts}^{-1} &=a_{ts}^{-1}a_{ts} =1 \ \ {\rm for} \ \ 1\leq s<t\leq n,\\
a_{ts}b_{rq}&=b_{rq}a_{ts} \ \ {\rm for} \ \ (t-r)(t-q)(s-r)(s-q)>0,\\
a_{ts}b_{ts}&=b_{ts}a_{ts}  \ \ {\rm for} \ \
1\leq s<t\leq n , \\
a_{ts}b_{sr} &=b_{tr}a_{ts}  \ \ {\rm for} \ \
1\leq r<s<t\leq n , \\
a_{sr}b_{tr} &=b_{ts}a_{sr}  \ \ {\rm for} \ \
1\leq r<s<t\leq n , \\
a_{tr}b_{ts}&=b_{sr}a_{tr}  \ \ {\rm for} \ \
1\leq r<s<t\leq n, \\
b_{ts}b_{rq}&=b_{rq}b_{ts} \ \ {\rm for} \ \ (t-r)(t-q)(s-r)(s-q)>0.
\end{cases}\label{eq:srebkl}
\end{equation}
\end{Theorem}

Now we consider the {\it positive singular braid monoid}   \BKL \ with respect 
to generators $a_{ts}$ and $b_{t,s}$ for $1\leq s < t \leq n$. Its
relations are (\ref{eq:srebkl}) except the one concerning the
invertibility of $a_{ts}$. Two positive words $A$ and $B$ in the
alphabet $a_{ts}$ and $b_{t,s}$
will be said to be {\it positively equivalent} if they are equal as 
elements of this monoid. In this case we shall write $A\doteq B$.

The {\it fundamental word} $\delta$ of Birman, Ko and Lee is given by
the formula
$$\delta \equiv a_{n(n-1)}a_{(n-1)(n-2)}
\dots a_{21} \equiv \sigma_{n-1} \sigma_{n-2} \dots \sigma_2\sigma_1.$$
Its divisibility by any generator $a_{ts}$, proved in \cite{BKL},
is convenient for us to be expressed in the following form.
\begin{Proposition} The fundamental word $\delta$ is positively equivalent to
 a word that begins or ends with any given generator $a_{ts}$.
The explicit expression for left divisibility is
\begin{equation*}
\delta \doteq a_{ts}a_{n(n-1)}a_{(n-1)(n-2)}
\dots a_{(t+1)s} a_{t(t-1)}\dots a_{(s+2)(s+1)} a_{s(s-1)}
\dots a_{21}.
\end{equation*}
\end{Proposition}
\begin{Proposition} For the fundamental word $\delta$ there are
the following formulae of commutation
\begin{equation*}
\begin{cases}
a_{ts} \delta &\doteq \delta a_{(t+1)(s+1)} \ \ \text{for}
\ \ 1\leq  s < t < n, \\
a_{ns} \delta &\doteq \delta a_{(s+1)1}, \\
b_{ts} \delta &\doteq \delta b_{(t+1)(s+1)} \ \ \text{for} \ \
 1\leq  s < t < n, \\
b_{ns} \delta &\doteq \delta b_{(s+1)1}.
\end{cases}
\end{equation*}
\end{Proposition}

\begin{figure}
\epsfbox{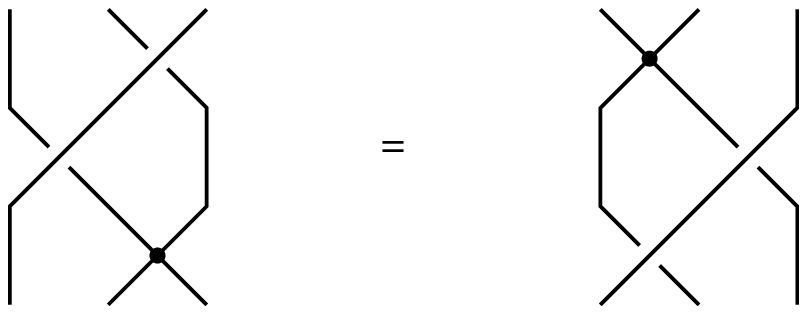}
\caption{}
\label{fi:bdelta}
\end{figure}

\begin{figure}
\epsfbox{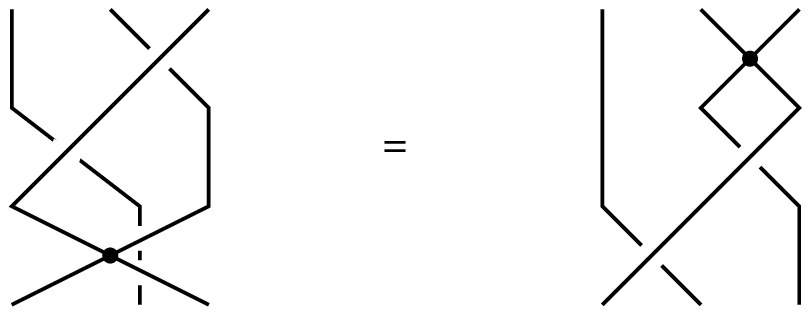}
\caption{}
\label{fi:bdeltam}
\end{figure}

Geometrically this commutation is shown on Figure~\ref{fi:bdelta}
and on Figure~\ref{fi:bdeltam}.

The analogues of the other results proved by Birman, Ko and Lee
 remain valid for the singular braid monoid.
They are proved in the work of V.~V.~Chaynikov \cite{Ch}.
\section{The work of V.~V.~Chaynikov\label{sec:chai}}

\subsection{Cancellation property}
Let $W_1,W_2 \in $\BKL. By a common multiple of $W_1,W_2$ (if it
exists) we mean a positive word $V \doteq W_1 V_1 \doteq W_2
  V_2$. 
  
Let
\begin{equation*} 
\de_{s_k...s_1} \equiv a_{s_k s_{(k-1)}}a_{s_{(k-1)}s_{(k-2)}} 
\dots a_{s_{21}},
\end{equation*}
 where $n \geq  s_k>s_{k-1}> \dots >s_{1} \geq 1$.
The word $\de_{s_k...s_1}$ is the least
   common multiple (l.c.m.) of the generators
$a_{ij}$, where $i,j \in \{ s_k, s_{(k-1)}, \dots s_1 \}$, see
\cite{BKL}, $\de \equiv \de _{n (n-1) \dots
1}$.

We denote the least common multiple of $X, Y$ by $X \vee Y$.
Define $(X \vee Y)_X^*$ and $(X \vee Y)_Y^*$ by equations:
$$ X \vee Y \doteq X (X \vee Y)_X^* \doteq Y(X \vee Y)_Y^*.$$
Similarly, we denote the greatest common divisor (g.c.d.) of $X , Y$
by $X \wedge Y$. The semigroup
$BKL_n^+$  is a lattice relative to $\vee, \wedge$ \cite{BKL}.

\begin{Remark}

We give the table of l.c.m. for  some pairs of generators below.
There does not exist $X \vee Y$ for the
rest pairs of  \BKL \ generators.

\begin{center}
\small
\begin{tabular}{|c|c|c|c|}\hline
  $X$ & $Y$ &$X \vee Y$& \\ \hline

$\hat{a}_{ts}$ & $\hat{a}_{rq}$ &$\hat{a}_{ts}( \hat{a}_{rq})\doteq 
\hat{a}_{rq} (\hat{a}_{ts})$ & $(t-r)(t-q)(s-r)(s-q)>0$ \\
\hline

& &$a_{ts}  (a_{sr}) \doteq a_{tr}  (a_{ts}) \doteq a_{sr}  (a_{tr})$& $t>s>r$ \\
\hline

$a_{ts}$ & $a_{rq}$ &$a_{ts}a_{tr}a_{sq} \doteq
a_{rq}a_{tq}(a_{rs}) \doteq \de_{trsq}$& $t>r>s>q$
\\ \hline

$a_{ts}$ & $b_{ts}$ &$a_{ts}  (b_{ts}) \doteq b_{ts} (a_{ts})$&
$t>s>r$
\\ \hline

$a_{sr}$ & $b_{tr}$ &$a_{ts}  (b_{sr}) \doteq b_{tr}(a_{ts})$& $t>s>r$ \\
\hline

$a_{sr}$ & $b_{ts}$ &$a_{sr} ( b_{tr}) \doteq b_{ts} (a_{sr})$& $t>s>r$ \\
\hline

$a_{ts}$ & $b_{sr}$ &$a_{ts} (a_{sr} b_{ts}) \doteq b_{sr} (\de_{tsr})$& $t>s>r$ \\
\hline

$a_{sr}$ & $b_{tr}$ &$a_{sr} ( a_{tr} b_{sr} ) \doteq b_{tr}( \de_{tsr})$& $t>s>r$ \\
\hline

$a_{tr}$ & $b_{ts}$ &$a_{tr} (a_{ts} b_{tr}) \doteq b_{ts} ( \de_{tsr})$& $t>s>r$ \\
\hline

$a_{ts}$ & $b_{rq}$ &$a_{ts} (\de_{trq} b_{ts}) \doteq b_{ts} (\de_{trsq})$& $t>r>s>q$ \\
\hline

$a_{rq}$ & $b_{ts}$ &$a_{rq} (a_{tq} a_{rs} b_{rq}) \doteq b_{ts} (\de_{trsq})$& $t>r>s>q$ \\
\hline

\end{tabular}
\end{center}

\smallskip\noindent
Here the symbol
  $\hat{a}_{ij} \in \{a_{ij},\ b_{ij}\}$,
 mean the same symbol in both parts of one equality.
\end{Remark}

 We call the pairs of generator from the table above
{\it admissible} and all other pairs {\it inadmissible}.
 Observe that pairs $\{ a_{ij}$, $a_{pm} \}$, $\{ a_{ij},b_{pm} \}$
are admissible  and $\{ b_{ij},b_{pm} \}$ is admissible  if and
only if $ b_{ij}b_{pm}=b_{pm}b_{ij}$ is the defining relation of
$SB_n$.
\begin{Theorem}[Left cancellation]\label{thm:L}
 i)Let $\{ x,y \}$ be an admissible  pair and $ xX \doteq yY$. Then
there exists a positive word $Z$ such that  $ xX \doteq yY \doteq
(x \vee y) Z$, where $ X \doteq (x \vee y)_x^* Z$ and $ Y \doteq
(x \vee y)_y^* Z$.

ii) If the pair $\{ x,y \}$ is inadmissible  then the equality $
xX \doteq yY$ is impossible (so there does not exist a common
multiple for  $\{ x,y \}$).
\end{Theorem}
Similarly we can obtain the Right cancellation property.
\begin{Corollary}
 If $A\doteq P$, $B\doteq Q$, $AXB \doteq PYQ$, then the equality $X\doteq Y$ holds in  \BKL.
\end{Corollary}
\begin{Corollary}Suppose that \De is the l.c.m. of the set of generators  
$\{ a_{i_1j_1}, \dots ,a_{i_pj_p}\}$
and $W$ is a positive word such that either
\begin{equation*}
W\doteq a_{i_1j_1} A_1\doteq a_{i_2j_2}
A_2\doteq\dots\doteq a_{i_pj_p}A_{p},
\end{equation*}
or
\begin{equation*}
W\doteq B_1a_{i_1j_1} \doteq B_2a_{i_2 j_2} \doteq\dots\doteq
B_p a_{i_p j_p},
\end{equation*}
then $W\doteq \de Z$ for some positive word $Z$.
\end{Corollary}

\begin{Corollary}[Embedding theorem]\label{thm:embedding}
The canonical homomorphism
$$SBKL_n^+ \to SB_n$$
is injective.
\end{Corollary}

\subsection{Word and conjugacy problems in $SB_n$\label{wcs}}
The word problem in $SB_n$ (in classical generators) was solved 
by R.~Corran \cite{Cor}, see also \cite{Ve10}.  
Let us fix an arbitrary linear order on the set of
generators of \BKL \  and extend it to the deg--lex order on
words of the generators of  \BKL. 
 With this order, we first order wwords  by total degree (the length of the word on given generators) and we break ties by the lex order. 
 By the  {\em base} of the positive word $W$ we mean the
 least
 (relative to the deg--lex order on the words on the generators of
\BKL)  word which represents the same element as $W$ in \BKL. 
Observe that this word is unique.
  If the positive word $A$ is not divisible by \De we denote its base by 
$\overline{A}$.

\begin{Theorem} Every word $W$ in $SB_n$ has a unique representation
of shape
 $ \de ^m \overline{A}$, where $m$ is an integer and $A$ is not divisible by \De.
\label{theo:garnf}
\end{Theorem}
This gives  a normal form for $SB_n$ in \bkl \  generators.
The process of computation of this normal form is the same as given by 
Garside \cite{Gar}. 
First, suppose that $P$ is any positive word in the generators \BKL. 
Among all positive
words positively equivalent to $P$ choose a word in the form
$\de^t A$ with $t$ maximal. Then $A$ is prime to $\de$ and we have
\begin{equation*}
P\doteq \de^t \overline{A}.
\end{equation*}
Now, let $W$ be an arbitrary word in $SB_n$. Then we may put
\begin{equation*}
W\equiv W_1(c_1)^{-1}W_2(c_2)^{-1} \dots (c_k)^{-1}W_{k+1},
\end{equation*}
where each $W_j$ is a positive word of length $\geq 0$, and $c_l$
are generators $a_{t,s}$, the only possible invertible generators. 
For each $c_l$ there exists a positive 
word $D_l$ such that $c_l D_l\doteq \de$, so that
$(c_l)^{-1} = D_l \de^{-1}$, and hence
\begin{equation*}
W = W_1 D_1 \de^{-1}W_2 D_2 \de^{-1} \dots
W_k D_k \de^{-1}W_{k+1}.
\end{equation*}
Moving the factors $\de^{-1}$ to the left, we obtain
$W=\de^k P$, where $P$ is positive, so we can express it in the
form $\de^t\overline{A}$ and finally we obtain the normal form
\begin{equation*}
W = \de^m \overline{A}.
\label{eq:gnf}
\end{equation*}

Let us consider the conjugacy problem.
We say that two elements  $u,v \in SB_n$ are conjugated if there
exists $g \in B_n$ such that $g^{-1}ug=v$. We denote this by $u \sim v$.

Let  $u$ be a positive word. Define the set of all positive
elements conjugated with $u$ as follows $C^+(u)=\{ v| v \sim u, v
\in SBKL_n^+ \}$.

The following properties are obvious and very close to the ones
proved in \cite{EM}, \cite{BKL}:

{\it i) The set of all positive words of limited length is finite.

ii)The set $C^+(u)$ is finite.

iii)The element $ \de^n $ generates the center of $SB_n$.}

Now fix two words $u,v \in SB_n$. We can
assume that they are positive (otherwise we multiply them by the
element  $ \de^{nk} $, where $k$ is big enough to cancel all
negative letters).

\begin{Theorem} The elements $u,v $ are conjugated if and only if the
sets $C^+(u)$ and $C^+(v)$ contain the same elements.
\end{Theorem}

There exists the following algorithm for constructing $C^+(u)$.
  Define  $C_0 ^+(u):=\{u\}$. 
If the set $C_{i}^+(u)$ is already constructed define 
$$C_{i+1} ^+(u):=\{ v^g \ | \ g \text{ divides} \ \de;\  v \in C_{i}^+
\} \cap SBKL_n^+.$$
The set $C_{k} ^+(u)$ stabilizes on the finite
step, so we put
$$C ^+(u):=\bigcup_{k\geq 0} C_k^+ (u).$$
\section{Inverse monoids\label{sec:invmon}}

The notion of {\it inverse semigroup} was introduced by V.~V.~Wagner in 1952 \cite{Wag}.
By definition it means that for any element $a$ of a semigroup (monoid)
$M$ there exists a unique element $b$ (which is called {\it inverse}) with the following two conditions: 
\begin{equation}
a = aba
\label{eq:reg_v_n}
\end{equation}
\begin{equation}
b = bab.
\label{eq:inv}
\end{equation}
Roots of this notion can be seen in the von Neumann regular rings \cite{v_N} 
where only one condition (\ref{eq:reg_v_n}) holds for non-necessary unique $b$,
or in the Moore-Penrose pseudoinverse for matrices \cite{Mo}, \cite{Pen}
where both conditions
(\ref{eq:reg_v_n}) and (\ref{eq:inv}) hold (and certain supplementary 
conditions also). See the books \cite{Pet} and \cite{Law} as general 
references for inverse semigroups.

The typical example of an inverse monoid is a monoid of partial (defined on a subset)
injections of a set. For a finite set this gives us the notion of a
{\it symmetric inverse monoid } $I_n$ which generalizes and includes the classical
symmetric group $\Sigma_n$. A presentation of symmetric inverse monoid was
obtained by L.~M.~Popova \cite{Po}, see also formulae 
(\ref{eq:invbrelations}\,-\ref{eq:syminvrelations})
below.

\begin{figure}
\epsfbox{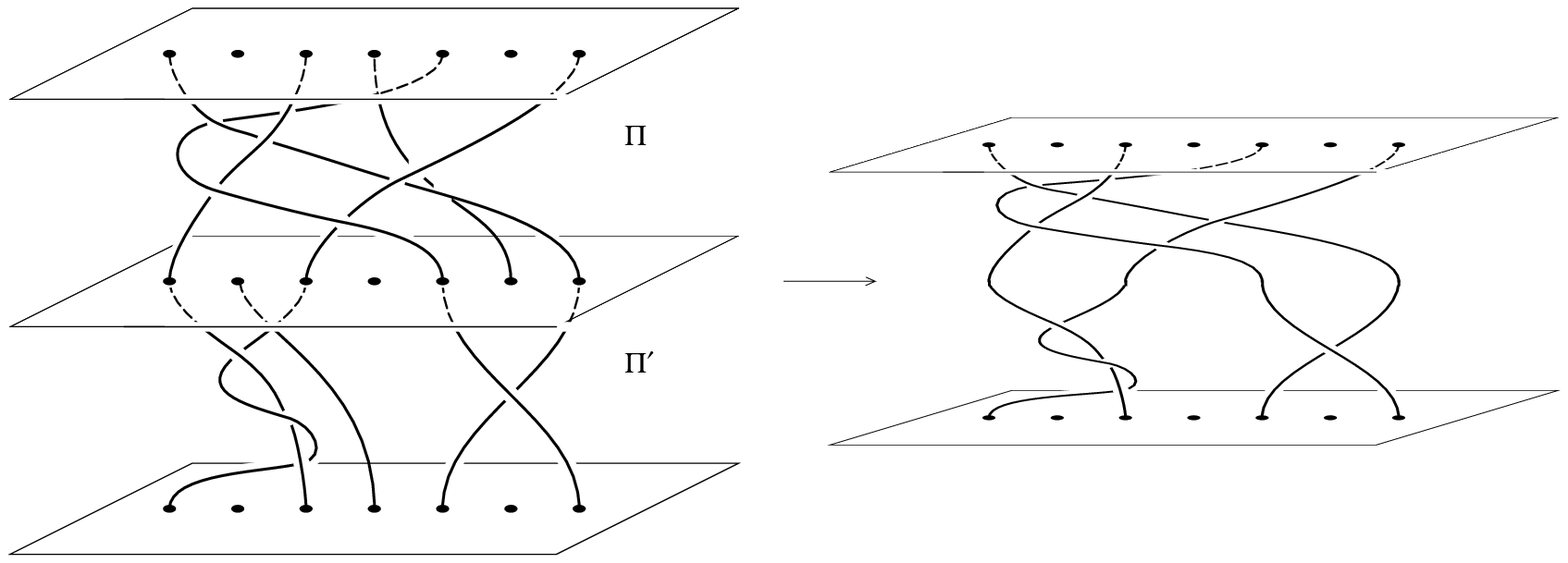}
\caption{} \label{fi:vv2}
\end{figure}

Recently the {\it inverse braid monoid }  $IB_n$ was constructed by 
D.~Easdown and T.~G.~Lavers \cite{EL}. 
It arises from 
a very natural operation on braids: deleting one or several strands.
By the application of this procedure to braids in
$Br_n$ we get {\it partial} braids  \cite{EL}.
The multiplication of partial braids is shown at Figure~\ref{fi:vv2}
At the last stage it is necessary to remove any arc that does not join the upper 
or lower planes.
 The set of all
partial braids  with this operation forms an inverse braid monoid $IB_n$. 

One of the motivations for studying $IB_n$ is that it is a natural setting for
the {\it Brunnian} (or {\it Makanin}) braids, which were also called  {\it smooth} braids by 
G.~S.~Makanin who first mentioned them in
\cite{Kou}, (page 78, question 6.23), and D.~L.~Johnson \cite{Joh1}. 
By the  usual definition a braid is   Brunnian if it becomes trivial
after deleting any strand, see formulae (\ref{eq:imak} - \ref{eq:mak}).
 According to the work of Fred Cohen, Jon Berrick, 
Wu Jie, Yang Loi Wong \cite{BCWW} Brunnian braids are connected with
homotopy groups of spheres.

The following presentation for the inverse braid monoid was obtained in
\cite{EL}. It has the generators $\sigma_i, \sigma_i^{-1} $, $i=1,\dots,n-1,$
$\epsilon$, which satisfy the braid relations (\ref{eq:brelations}
and the following relations:
\begin{equation}
 \begin{cases} 
&\sigma_i\sigma_i^{-1}=\sigma_i^{-1}\sigma_i =1, \ \text {for \ all} \ i, \\
&\epsilon \sigma_i  =\, \sigma_i \epsilon \ \ \text {for } i\geq 2,   \\ 
&\epsilon\sigma_1 \epsilon  = \sigma_{1} \epsilon \sigma_1 \epsilon = 
\epsilon\sigma_{1} \epsilon \sigma_1, \\
&\epsilon = \epsilon^2 = \epsilon \sigma_1^2= \sigma_1^2 \epsilon.
\end{cases} \label{eq:invbrelations}
\end{equation}
Geometrically the generator $\epsilon$ means that the first strand in 
the trivial braid is absent.
 
If we replace the first relation in (\ref{eq:invbrelations})
 by the following set of relations
\begin{equation}
\sigma_i^2 =1, \ \text {for \ all} \ i, \\
 \label{eq:syminvrelations}
\end{equation}
and delete the
superfluous relations 
$$\epsilon =  \epsilon \sigma_1^2= \sigma_1^2 \epsilon, $$
we get a presentation of the symmetric inverse monoid $I_n$ \cite{Po}. 
We also can simply add the relations (\ref{eq:syminvrelations})
if we do not worry about redundant relations.
We get a canonical map \cite{EL}
\begin{equation}
\tau_n: IB_n\to I_n
 \label{eq:tauIBn}
\end{equation}
which is a natural extension of the corresponding map for the braid 
and symmetric groups.

More balanced relations for the inverse braid monoid were 
obtained in \cite{Gil}.
Let $\epsilon_i$ denote the braid which is obtained from the trivial 
by deleting of the $i$th strand,
formally:
\begin{equation*}
\begin{cases} \epsilon_1 &= \epsilon, \\
\epsilon_{i+1} &= \sigma_i^{\pm 1}\epsilon_{i}\sigma_i^{\pm 1}.  
\end{cases} \end{equation*} 
So, the generators are: $\sigma_i, \sigma_i^{-1} $, $i=1,\dots,n-1,$
$\epsilon_i$, $i=1,\dots,n$, and relations  are the following:

\begin{equation}
 \begin{cases} 
&\sigma_i\sigma_i^{-1}=\sigma_i^{-1}\sigma_i =1, \ \text {for all} \ i, \\
&\epsilon_j \sigma_i =\, \sigma_i \epsilon_j \ \ \text {for } \ j \not= i, i+1, \\ 
&\epsilon_i\sigma_i =  \sigma_{i} \epsilon_{i+1},  \\
&\epsilon_{i+1}\sigma_i =  \sigma_{i} \epsilon_{i},  \\
&\epsilon_i = \epsilon_i^2 , \\
& \epsilon_{i+1} \sigma_i^2= \sigma_i^2 \epsilon_{i+1} = \epsilon_{i+1}, \\
&\epsilon_i \epsilon_{i+1} \sigma_i = \sigma_{i} \epsilon_i \epsilon_{i+1}
=\epsilon_i\epsilon_{i+1},
\end{cases} \label{eq:invbrelations2}
\end{equation}
plus the braid relations (\ref{eq:brelations}).

\subsection{Inverse reflection monoid of type $B$\label{sec:typeb}}

It can be defined in the same way as the corresponding Coxeter group (\ref{weylb})
as the {\it monoid of partial signed permutations} $I(B_n)$: 
\begin{multline*}
I(B_n)=\{\sigma \text{ is a partial 
bijection of } SN: (-x)\sigma  =-(x)\sigma \text{ for} \ x\in SN \\
\text {and } x\in \operatorname{dom} \sigma \ \text {if and only if}
-x\in \operatorname{dom} \sigma\}, 
\end{multline*}
where $\operatorname{dom} \sigma$ means  domain of definition of
the monomorphism $\sigma$. This monoid was studied in \cite{EF}.

\section{Properties of inverse braid monoid\label{sec:prop}}

In relations  (\ref{eq:invbrelations}) we have one generator
for the idempotent part and $n-1$ generators for the group part.
If we minimize the number of generators of the group part and take the 
presentation  (\ref{eq:2relations}) for the braid group we get a presentation
of the inverse braid monoid with generators $\sigma_1, \sigma$,  $\epsilon$,
and relations:
\begin{equation*}
 \begin{cases} 
&\sigma_1\sigma_1^{-1}=\sigma_1^{-1}\sigma_1 =1,  \\
&\sigma\sigma^{-1}=\sigma^{-1}\sigma =1,  \\
&\epsilon \sigma^{i}\sigma_1\sigma^{-i}  
=\, \sigma^{i}\sigma_1\sigma^{-i} \epsilon \ \ \text {for } 1\leq i \leq n-2,   \\ 
&\epsilon\sigma_1 \epsilon  = \sigma_{1} \epsilon \sigma_1 \epsilon = 
\epsilon\sigma_{1} \epsilon \sigma_1, \\
&\epsilon = \epsilon^2 = \epsilon \sigma_1^2= \sigma_1^2 \epsilon,
\end{cases}
\end{equation*}
plus (\ref{eq:2relations}).

Let $\Gamma$ be a normal planar graph (see Introduction). Let us add new generators $\epsilon_v$ which 
correspond to each vertex of the graph $\Gamma$. Geometrically it means the
absence in the trivial braid of one strand corresponding to the vertex $v$.
We orient the graph $\Gamma$ arbitrarily and so we get a starting $v_0=v_0(e)$
and a terminal $v_1=v_1(e)$ vertex for each edge $e$.
Consider the following relations
\begin{equation}
 \begin{cases} 
&\sigma_e\sigma_e^{-1}=\sigma_e^{-1}\sigma_e =1, \ \text {for  all  edges  of} \  
\Gamma, \\
&\epsilon_v \sigma_e  =\, \sigma_e \epsilon_v, \ \ \text {if  the vertex } \ 
v  \ \text{and  the  edge} \ e \ \text{do  not  intersect}, \\ 
&\epsilon_{v_0}\sigma_e =  \sigma_{e} \epsilon_{v_1}, \ \text{where} \
v_0=v_0(e), \ v_1= v_1(e), \\
&\epsilon_{v_1}\sigma_e =  \sigma_{e} \epsilon_{v_0},  \\
&\epsilon_v = \epsilon_\nu^2 , \\
& \epsilon_{v_i} \sigma_e^2= \sigma_e^2 \epsilon_{v_i} = \epsilon_{v_i}, 
\ \ i=0,1, \\
&\epsilon_{v_0}\epsilon_{v_1} \sigma_e = \sigma_{e} \epsilon_{v_0} \epsilon_{v_1}
=\epsilon_{v_0}\epsilon_{v_1}.
\end{cases} \label{eq:invbserel}
\end{equation}

\begin{Theorem} We get a Sergiescu graph presentation of the inverse braid monoid 
$IB_n$ if we add to the graph presentation of the braid group $Br_n$ 
 relations (\ref{eq:invbserel}).
\end{Theorem} 
Let $E F_n$ be a monoid of partial isomorphisms of a free group
$F_n$ defined as follows. Let $a$ be an element of the symmetric
inverse monoid $I_n$, $a\in I_n$, $J_k =\{j_1, \dots, j_k\}$ 
is the image of $a$, and elements $i_1, \dots, i_k$ belong to
domain of the definition of $a$. The monoid $E F_n$
consists of isomorphisms of free subgroups
$$<x_{i_1}, \dots, x_{i_k}> \, \to \ <x_{j_1}, \dots, x_{j_k}>$$
such that
$$f_a :x_i\mapsto w_i^{-1} x_{a(i)}w_i, $$
if $i$ is among $i_1, \dots, i_k$ and not defined otherwise and 
$w_i$ is a word on $x_{j_1}, \dots, x_{j_k}$.
The composition of $f_a$ and $g_b$, $a, b\in I_n$, 
is defined for $x_i$ belonging to the domain of $a\circ b$.
We put $x_{j_m}=1$ in a word $w_i$ if $x_{j_m}$ does not belong
to the domain of definition of $g$.
 We define a map $\phi_n$ from $IB_n$ to $E F_n$
expanding the canonical inclusion 
\begin{equation*}
Br_n \to  \operatorname{Aut} F_n
 \end{equation*}
by the condition that $\phi_n(\epsilon)$ 
as a partial isomorphism of $F_n$ is given by the
formula
\begin{equation} 
\phi_n(\epsilon)(x_i) = \begin{cases}
x_i {\text{ if} } \ i\geq 2 , \\ 
{\text {not defined,  if }} i=1 . 
\end{cases} \label{eq:endf1}
\end{equation}

Using the presentation (\ref{eq:invbrelations}) we see that $\phi_n$ is 
correctly defined homomorphism of monoids
\begin{equation*}
\phi_n: IB_n \to  E F_n.
 \end{equation*}
 
\begin{Theorem} The homomorphism $\phi_n$ is a monomorphism.
\label{theo:isoendo} 
\end{Theorem} 
Theorem~\ref{theo:isoendo}  gives also a possibility to interpret the inverse 
braid monoid as
a monoid of isotopy classes of maps. As usual consider a disc $D^2$ with $n$ 
fixed points. Denote the set of these points by $Q_n$.
The fundamental group of $D^2$ with these points deleted is isomorphic to $F_n$. 
Consider homeomorphisms
 of $D^2$ onto a copy of the same disc with the condition that
only $k$ points of $Q_n$,  $k \leq n$ (say $i_1, \dots, i_k$) are mapped
bijectively onto the $k$ points (say $j_1, \dots, j_k$) of the second copy 
of $D^2$. 
Consider the isotopy 
classes of such homeomorphisms and denote such set  by $IM_n(D^2)$. Evidently 
it is a monoid.
\begin{Theorem} The monoids $IB_n$ and  $IM_n(D^2)$
are isomorphic.
\end{Theorem} 
These considerations can be generalized to the following definition.
Consider a surface $S_{g,b,n}$ of the genus $g$, $b$ boundary components and
with a chosen set $Q_n$ of $n$ fixed interior points. Let $f$ be a homeomorphism of
 $S_{g,b,n}$ 
which maps  $k$ points,  $k \leq n$, from $Q_n$:  $\{i_1, \dots, i_k\}$ to 
$k$ points 
$\{j_1, \dots, j_k\}$ also from $Q_n$. In the same way let $h$  
be a homeomorphism of $S_{g,b,n}$ 
which maps  $l$ points, $l\leq n$, from $Q_n$, say $\{s_1, \dots, s_l\}$ to $l$ points 
$\{t_1, \dots, t_l\}$ again from $Q_n$. Consider the intersection of the sets
$\{j_1, \dots, j_k\}$ and $\{s_1, \dots, s_l\}$,  let it be the set of cardinality $m$,
it may be empty. Then the composition of $f$ and $h$ maps $m$ points of $Q_n$
to $m$ points (may be different) of $Q_n$. If $m=0$ then the composition does not 
take into account the set $Q_n$. Denote the set of isotopy classes of such maps
by $\mathcal I \mathcal {M}_{g,b,n}$. This standard composition of $f$ and $g$
as maps  defines a structure of monoid
on $\mathcal I \mathcal {M}_{g,b,n}$.
\begin{Proposition} The monoid
$\mathcal I \mathcal {M}_{g,b,n}$ is inverse.
\end{Proposition}
We call the monoid $\mathcal I \mathcal {M}_{g,b,n}$ the {\it inverse mapping class monoid}.
If $g=0 $ and $b=1$ we get the inverse braid monoid. In 
the general case $\mathcal I \mathcal {M}_{g,b,n}$ the role of the empty braid plays the
mapping class group  $ \mathcal {M}_{g,b}$ (without fixed points). 

We remind that a monoid $M$ is {\it factorisable} if $M= EG$ where $E $ 
is a set of 
idempotents of $M$ and $G$ is a subgroup of $M$. 

\begin{Proposition} The monoid
$\mathcal I \mathcal {M}_{g,b,n}$ 
can be written in the form
$$\mathcal I \mathcal {M}_{g,b,n} = E \mathcal {M}_{g,b,n},$$
where $E $ is a set of 
idempotents of $\mathcal I \mathcal {M}_{g,b,n}$  and $\mathcal {M}_{g,b,n}$ is 
the corresponding mapping class group. So this monoid is factorisable.
\end{Proposition}

Let $\Delta$ be the Garside's   fundamental word   in the braid 
group $Br_{n}$ \cite{Gar}. It can be defined by the formula:
$$\Delta = \sigma_1 \dots \sigma_{n-1} \sigma_1 \dots \sigma_{n-2} \dots  
\sigma_1 \sigma_2 \sigma_1.$$
\begin{Proposition} The generators $\epsilon_i$ commute with $\Delta$ in
the following way:
\begin{equation*}
\epsilon_i\Delta = \Delta \epsilon_{n+1-i}.
\end{equation*}
\end{Proposition}
\begin{Proposition} The center of $IB_n$ consists of the union of the center of the
braid group $Br_n$ (generated by $\Delta^2$) and the empty braid 
$\varnothing = \epsilon_1 \dots \epsilon_n$.
\end{Proposition}

Let $\mathcal E$ be the monoid generated by one idempotent generator 
$\epsilon$ . 
\begin{Proposition} The abelianization of  $IB_n$ is isomorphic to
an abelian monoid ${AB}$ generated (as an abelian monoid)
by elements $\epsilon$, $\alpha$  and $-\alpha$, subject to the following relations
\begin{equation*}
\begin{cases}
\alpha+(-\alpha)=0,\\
2\epsilon =\epsilon, \\
\epsilon + \alpha =\epsilon. 
\end{cases}
\end{equation*}
So, it is isomorphic to the quotient-monoid
of $\mathcal E \oplus \Z$ by the relation 
$\epsilon +  1 =\epsilon$.
The canonical map of abelianization
\begin{equation*}
a: {IB}_{n} \to {AB}
\end{equation*}
is given by the formula:
\begin{equation*}
\begin{cases}
a(\epsilon_i) = \epsilon ,\\ 
a(\sigma_i) = { \alpha} .
\end{cases}
\end{equation*}
\end{Proposition}

Let $\epsilon_{k+1, n}$ denote the partial braid with the trivial first $k$ strands 
and the absent rest $n-k$ strands. It can be expressed using the generator $\epsilon$ 
or the generators $\epsilon_i$ as follows
\begin{equation} \epsilon_{k+1, n}=
\epsilon\sigma_{n-1}\dots\sigma_{k+1}\epsilon \sigma_{n-1}\dots\sigma_{k+2}
\epsilon\dots \epsilon\sigma_{n-1}\sigma_{n-2}\epsilon \sigma_{n-1}\epsilon,
\end{equation}
\begin{equation} \epsilon_{k+1, n}=
\epsilon_{k+1}\epsilon_{k+2} \dots \epsilon_{n},
 \label{eq:espki}
\end{equation}
It was proved in \cite{EL} the 
every partial braid has a representative of the form
\begin{equation} \sigma_{i_1}\dots\sigma_{1}\dots \sigma_{i_k}\dots\sigma_{k}
\epsilon_{k+1, n}x \epsilon_{k+1, n}\sigma_{k}\dots\sigma_{j_k}\dots\sigma_{1}
\dots\sigma_{j_1},
 \\ 
\label{eq:form_inv}\end{equation}
\begin{equation} k\in \{0,\dots, n\}, x\in Br_k,
0\leq i_1<\dots<i_k\leq n-1  \ 
\text{and} \  0\leq j_1<\dots<j_k\leq n-1.
\end{equation}
Note that in the formula (\ref{eq:form_inv}) we can  delete one of the 
$\epsilon_{k+1,n}$, but we shall use the form (\ref{eq:form_inv}) because of 
convenience: two symbols $\epsilon_{k+1,n}$ serve as markers to distinguish
the elements of $Br_k$.
We can put the element $x\in Br_k$ in the Markov normal form  \cite{Mar2}
and get the corresponding {\it Markov normal form for the inverse braid monoid} 
$IB_n$. 

Among positive words on the alphabet  $\{\sigma_1 \dots \sigma_n\}$ let us 
introduce a lexicographical ordering with the condition that  
$\sigma_1 < \sigma_2 < \dots < \sigma_n $. For a positive word $V$ the 
\emph{base} of $V$ is the smallest positive word which is positively equal
to $V$. The base is uniquely determined. If a positive word $V$ is prime 
to $\Delta$, then for the base of $V$ the notation $\overline{V}$  will
be used (compare with subsection~\ref{wcs}).
\begin{Theorem}
 Every word $W$ in $IBr_{n}$ can be uniquely written in
the form 
\begin{equation} \sigma_{i_1}\dots\sigma_{1}\dots \sigma_{i_k}\dots\sigma_{k}
\epsilon_{k+1, n}x \epsilon_{k+1, n}\sigma_{k}\dots\sigma_{j_k}\dots\sigma_{1}
\dots\sigma_{j_1},
 \\ 
\end{equation}
\begin{equation} k\in \{0,\dots, n\}, x\in Br_k,
0\leq i_1<\dots<i_k\leq n-1  \ 
\text{and} \  0\leq j_1<\dots<j_k\leq n-1.
\end{equation} 
where $x$ is written in the Garside normal form for $Br_k$
$$\Delta^m \overline{V},$$
where $m$ is an  integer.
\label{Theorem:garnfm}
\end{Theorem}
Theorem~\ref{Theorem:garnfm} is evidently true also for the presentation with 
$\epsilon_i$,
$i= 1,\dots n$. In this case the elements $\epsilon_{k+1,n}$ are expressed
by (\ref{eq:espki}). 

We call the form of a word $W$ established in Theorem~\ref{Theorem:garnfm} the
\emph{Garside left normal form for the inverse braid monoid} $IB_n$
and the 
index $m$ we call the \emph{power} of $W$. In the same way we can define the
\emph{Garside right normal form for the inverse braid monoid}  and 
the corresponding
variant of Theorem~\ref{Theorem:garnfm} is true. 
\begin{Theorem} 
The necessary and sufficient condition for two words in  
$IB_{n}$ to be equal is that their Garside normal forms are identical.
The Garside normal form  gives a solution to the word problem in the braid
group.
\label{Theorem:gws}
\end{Theorem}

Garside normal form for the braid groups was detailed in the subsequent 
works of S.~I.~Adyan \cite{Ad},
 W.~Thurston \cite{E_Th}, E.~El-Rifai and H.~R.~Morton \cite{EM}. 
Namely, there was introduced 
the \emph{left-greedy form} (in the terminology of
W.~Thurston \cite{E_Th}) 
\begin{equation*}
\Delta^t A_1 \dots A_k,
\end{equation*}
where $A_i$ are the successive
possible longest \emph{fragments of the word} $\Delta$ (in the terminology 
of S.~I.~Adyan \cite{Ad}) or \emph{positive permutation braids} (in the 
terminology of E.~El-Rifai and H.~R.~Morton \cite{EM}). In the same 
way one defines the \emph{right-greedy form} is defined. These greedy forms are
defined for the inverse braid monoid in the  same way.

Let us consider the elements $m\in IB_n$ satisfying  the equation:
\begin{equation}
\epsilon_i m = \epsilon_i.
\label{eq:imak}
\end{equation}
Geometrically this means that removing the strand  (if it exists) that 
starts at the point with the number $i$ we get a trivial braid on the 
remaining $n-1$ strands. It is equivalent to the condition
\begin{equation}
m \epsilon_{\tau (m)(i)}  = \epsilon_{\tau(m)(i)},
\label{eq:imak2}
\end{equation}
where $\tau$ is the canonical map to the symmetric monoid (\ref{eq:tauIBn}).
With the exception of $\epsilon_i$ itself all such elements belong to
$Br_n$. We call such braids as $i$\,-{\it Brunnian} and denote the subgroup
of $i$\,-Brunnian braids by $A_i$. The  subgroups $A_i$, $i=1, \dots, n$,
 are conjugate
\begin{equation} 
A_i = \sigma_{i-1}^{-1} \dots \sigma_1^{-1} A_1\sigma_1 \dots \sigma_{i-1}
\label{eq:Aimak}
\end{equation} 
free subgroups. The group $A_1$ is freely generated by the set 
$\{x_1, \dots, x_{n-1}\}$ \cite{Joh1}, where 
\begin{equation}
x_i = 
\sigma_{i-1}^{-1} \dots \sigma_1^{-1} \sigma_1^2\sigma_1 \dots \sigma_{i-1}.
\label{eq:freegmak}
\end{equation} 
The intersection of all subgroups of $i$\,-Brunnian braids is the group of 
Brunnian braids
\begin{equation} 
Brunn_n =\cap_ {i=1}^{n} A_i.
\label{eq:mak}
\end{equation}
That is the same as $m\in Brunn_n$ if and only if the equation (\ref{eq:imak})
holds for all $i$.

\section{Monoids of partial generalized braids\label{sec:m_gen_br}}

Construction of partial braids can be applied to various generalizations of 
braids, namely to those where geometric or diagrammatic construction of braids 
takes place.
Let $\Sigma_g$ be a surface of genus $g$ possibly with boundary components and 
punctures. We consider partial braids lying in a layer between two such surfaces:
$\Sigma_g\times I$ and take a set of isotopy classes of such braids.   
We get a monoid of partial braid on a surface $\Sigma_g$,  denote it by $IB_n(\Sigma_g)$. 
An interesting case is when the surface is a sphere $S^2$. So our partial 
braids are lying in a layer between two concentric spheres. 

\begin{Theorem} We get a presentation of the monoid $IB_n(S^2)$ if we add to
the presentation (\ref{eq:invbrelations}) or to the presentation (\ref{eq:invbrelations2})
of  $IB_n$ the sphere relation (\ref{eq:spherelation}).
It is a factorisable inverse monoid.
\end{Theorem} 
The monoid $IB(B_n)$ of partial braids of the type $B$
 can be considered also as a submonoid of $IB_{n+1}$ consisting of
partial braids with the first strand fixed. An interpretation  as
a monoid of isotopy classes of homeomorphisms is possible as well. Consider a disc $D^2$ with given $n+1$  points. Denote the set of these points by
 $Q_{n+1}$.
Consider homeomorphisms
 of the disc $D^2$ onto a copy of the same disc with the condition that
 the first point is always mapped into itself and among the other $n$
 points
only $k$ points,  $k \leq n$ (say $i_1, \dots, i_k$) are mapped
bijectively onto the $k$ points (say $j_1, \dots, j_k$) of the 
set $Q_{n+1}$ (without the first point) of second copy of the disc $D^2$. 
The isotopy classes of such homeomorphisms form the monoid $IB(B_n)$.
\begin{Theorem}  We get a presentation of the monoid $IB(B_n)$ if we add to
the presentation (\ref{eq:invbrelations}) or the presentation (\ref{eq:invbrelations2})
of  $IB_n$ one generator $\tau$, the type $B$ relation
(\ref{eq:relB}) and the following relations
\begin{equation}
 \begin{cases} 
&\tau\tau^{-1}=\tau^{-1}\tau =1,  \\
&\epsilon_1\tau = \tau\epsilon_1 = \epsilon_1.
\end{cases} \label{eq:IBB}
\end{equation}
It is a factorisable inverse monoid.
\label{con}
\end{Theorem} 
\begin{Remark} 
Theorem~\ref{con} can be naturally generalized for partial braids in 
handlebodies \cite{Ve1}.
\end{Remark} 
 We define an action of the monoid $ {IB}(B_n)$ on the set $SN$ (see 
 subsection~\ref{subsection:abbg}) by partial 
 isomorphisms  as follows
 \begin{equation}
 \sigma_i(\delta_j v_j) =\begin{cases} \delta_{i} v_{i+1}, \ \text{if } j=i,\\
  \delta_{i+1} v_i, \ \text{if } j=i+1,\\
  \delta_j v_j,  \ \text{if } j\not=i,i+1, \\
  \end{cases}
  \label{eq:sigmai}
 \end{equation}
 \begin{equation}
 \tau(\delta_j v_j) =\begin{cases} -\delta_{1} v_{1}, \ \text{if } j=1,\\
  \delta_j v_j,  \ \text{if } j\not=1,\\
  \end{cases}
 \end{equation} 
 \begin{equation}
 \operatorname{dom}\epsilon =\{\delta_2 v_2, \dots, \delta_n v_n\},
 \end{equation}  
 \begin{equation}
 \epsilon(\delta_j v_j) =  \delta_j v_j,  \ \text{if } j=2, \dots, n,\\
 \end{equation} 
  \begin{equation}
 \operatorname{dom}\epsilon_i =\{\delta_1 v_1, \dots, {\widehat{\delta_i v_i}}, 
 \dots, \delta_n v_n\},
 \end{equation}  
 \begin{equation}
 \epsilon_i(\delta_j v_j) =  \delta_j v_j,  \ \text{if } j=1, \dots, \widehat{i},
 \dots, n.\\
  \label{eq:epsiloni}
 \end{equation} 
 Direct checking shows that the relations of the inverse braid monoid of type 
 $B$ are satisfied by the corresponding compositions of partial isomorphisms defined by
  $\sigma_i$, $\tau$  and $\epsilon_i$.
\begin{Theorem}  The action given by the formulae (\ref{eq:sigmai} -  
\ref{eq:epsiloni}) defines a homomorphism of inverse monoids
$\rho_B: IB(B_n)\to I(B_n)$ such that the 
following diagram commutes
\begin{equation}
\CD
Br(B_n)@>>>W(B_n)\\
\downarrow&&\downarrow\\
IB(B_n)@>\rho_B>> I(B_n) 
\endCD
\label{eq:cod}
\end{equation}
(where the vertical arrows mean inclusion of the group of invertible elements
into a monoid).
\end{Theorem} 
\begin{Theorem} The homomorphism 
$\rho_B: IB(B_n)\to I(B_n)$
is an epimorphism.
We get a presentation of the monoid $ I(B_n)$ if in
the presentation of ${IB}(B_n)$ 
 we replace the first relation in (\ref{eq:invbrelations})
 by the following set of relations
\begin{equation*}
\sigma_i^2 =1, \ \text {for \ all} \ i, \\
\end{equation*}
and delete the superfluous relations 
$$\epsilon =  \epsilon \sigma_1^2= \sigma_1^2 \epsilon, $$
and 
we replace the first relation in (\ref{eq:IBB})
 by the following relation
\begin{equation*}
\tau^2 =1.  \\
\end{equation*}
\end{Theorem}

We remind that  $\mathcal E$ denotes the monoid generated by one 
idempotent generator 
$\epsilon$. 
\begin{Proposition} The abelianization $Ab(IB(B_n))$ of the monoid
 $IB(B_n)$ is isomorphic to the monoid
$\mathcal E \oplus \Z^2$, factorized by the relations
\begin{equation*}
\begin{cases}
\epsilon + \tau =\epsilon,\\
\epsilon + \sigma =\epsilon, \\
\end{cases}
\end{equation*}
where $\tau$ and $\sigma$ are generators of $\Z^2$. 
 The canonical map of abelianization
\begin{equation*}
a: IB(B_n) \to Ab(IB(B_n))
\end{equation*}
is given by the formulae:
\begin{equation*}
\begin{cases}
a(\epsilon_i) = \epsilon ,\\
a(\tau) = \tau ,\\
a(\sigma_i) = \sigma .
\end{cases}
\end{equation*}
The canonical map from $Ab(IB(B_n))$ to $Ab(I(B_n))$ consists of 
factorizing $\Z^2$ modulo $2$.
\end{Proposition}
 
Let $BP_n$ be the braid-permutation group (see subsection~\ref{bp}).  
Consider the image of monoid $I_n$ in $EF_n$
by the map  defined by the formulae (\ref{eq:perm}), (\ref{eq:endf1}).
Take also the  monoid $IB_n$  lying in $EF_n$
under the map $\phi_n$ of Theorem~(\ref{theo:isoendo}). We define the 
{\it braid-permutation} 
monoid as a submonoid of $EF_n$ generated by both
images of $IB_n$ and $I_n$ and denote it by $IBP_n$. It can be also defined
by the diagrams of partial welded braids. 
\begin{Theorem} We get a presentation of the monoid $IBP_n$ if we add to
the presentation of  $BP_n$ the generator $\epsilon$, relations 
(\ref{eq:invbrelations})
and the analogous relations between
$\xi_i$ and $\epsilon$, or 
generators $\epsilon_i$, $ 1\leq i \leq n$
relations (\ref{eq:invbrelations2}) and the analogous relations between
$\xi_i$ and $\epsilon_i$. 
It is a factorisable inverse monoid. Monoid $IBP_n$ is isomorphic to the monoid
$EF_n$ of partial isomorphisms of braid-conjugation type.
\end{Theorem} 

The virtual braids \cite{Ve8} can be defined by the plane diagrams with real 
and virtual crossings. The corresponding Reidemeister moves are the same as 
for the welded braids of the  braid-permutation group with one exception.
 The forbidden move
corresponds to the last mixed relation for the braid-permutation group
(\ref{eq:mixperm}).
This allows to define the partial virtual braids and the corresponding
monoid $IVB_n$. So the mixed relation for $IVB_n$ have the form:
\begin{equation} \begin{cases} \sigma_i \xi_j
&=\xi_j \sigma_i, \ \text {if} \  |i-j| >1,
\\ \xi_i \xi_{i+1} \sigma_i &= \sigma_{i+1} \xi_i \xi_{i+1}.
\end{cases} \label{eq:mixvir}
\end{equation}
\vglue0.01cm
\centerline {The mixed relations for virtual braids}
\smallskip

\begin{Theorem} We get a presentation of the monoid $IVB_n$ if we delete
the last mixed relation in the presentation of  $IBP_n$, that is 
replace the relations (\ref{eq:mixperm}) by (\ref{eq:mixvir})  
It is a factorisable inverse monoid. The canonical epimorphism 
$$ IVB_n\to IBP_n$$
is evidently defined.
\end{Theorem} 

The constructions of singular braid monoid $SB_n$ (see subsection~\ref{sbm})
are geometric, so we can easily get the analogous
monoid of partial singular braids $PSB_n$.
\begin{Theorem} We get a presentation of the monoid $PSB_n$ if we add to
the presentation of  $SB_n$ the  generators $\epsilon_i$, $ 1\leq i \leq n$,
relations (\ref{eq:invbrelations2}) and the analogous relations between
$x_i$ and $\epsilon_i$. 
\end{Theorem} 
\begin{Remark} The monoid $PSB_n$  is not neither factorisable nor  inverse.
\end{Remark}

The construction of braid groups on graphs \cite{Ghr}, \cite{FarS} is geometrical
so, in the same way as for the classical braid groups we can define {\emph {partial
braids on a graph }}$\Gamma$ and the {\emph {monoid of partial braids on a graph}}
$\Gamma$ which will be evidently inverse, so we call it as
{\emph {inverse braid monoid on the graph}} $\Gamma$ and we denote it as
$IB_n\Gamma$.


\begin{thebibliography}{References}
\bibitem{Ad}
\emph{S.~I.~Adyan}, Fragments of the word $\Delta $ in the braid group. 
(Russian) Mat. Zametki 36 (1984), no. 1, 25--34. 
\bibitem{Art1}
\emph{ E. Artin},
 Theorie der Z\"opfe. 
 Abh. Math. Semin. Univ. Hamburg, 1925,
v. 4, 47--72.
\bibitem{Bae}
\emph{J.~C.~Baez}, Link invariants of finite type and perturbation theory. 
Lett. Math. Phys. 26 (1992), no. 1, 43--51.
\bibitem{BN1}
\emph{D.~Bar-Natan}, Non-associative tangles.  Geometric topology 
(Athens, GA, 1993),  139--183, AMS/IP Stud. Adv. Math., 2.1, Amer. Math. Soc., 
Providence, RI, 1997.
\bibitem{BN2}
\emph{D.~Bar-Natan},
On associators and the Grothendieck-Teichmuller group. I.  Selecta Math. (N.S.)  
4  (1998),  no. 2, 183--212.
\bibitem{BaMi}
\emph{V.~Bardakov, R.~Mikhailov},
On certain questions of the free group automorphisms theory.
Comm. Algebra 36 (2008), no. 4, 1489--1499.
\bibitem{beltes} 
\emph{P.~Bellingeri}, 
Surface braid groups and polynomial link invariants, 
Thesis, Univ. Grenoble I (2003. 
\bibitem{BelV}
\emph{P.~Bellingeri, V.~Vershinin}, Presentations of surface  braid groups by graphs,
 Fund. Math. Vol. 188, December, 2005, 1-20.
\bibitem{BP}
\emph{ B.~Berceanu, S.~Papadima},
Universal representations of braid and braid-permutation groups.
J. Knot Theory Ramifications 18 (2009), no. 7, 999--1019.
\bibitem{BCWW}
\emph{J.~A.~Berrick, F.~R.~Cohen, Y.~L.~Wong and J.~Wu}, Configurations, b
raids, and homotopy groups,  J. Amer. Math. Soc.  19  (2006),  no. 2, 265--326. 
\bibitem{Be}
\emph{D.~Bessis}, Finite complex reflection arrangements are $K(\pi,1)$, arXiv:math.GT/0610777
\bibitem{BM}
\emph{D.~Bessis, J.~Michel}, Explicit presentations for exceptional braid groups. 
Experiment. Math. 13 (2004), no. 3, 257--266. 
\bibitem{Bir2}
\emph{J.~S.~Birman}, New points of view in knot theory,
Bull. Amer. Math. Soc.
 1993, 28{\rm , No 2}, 253--387.
\bibitem{BKL}
\emph{J.~S.~Birman,  K.~H.~Ko, S.~J.~Lee}, A new approach to the 
word and conjugacy problems in the
   braid groups. Adv. Math. 139 (1998), no. 2, 322--353. 
\bibitem{Bo}
\emph{N.~Bourbaki}, Groupes et alg\`ebres de Lie, Chaps. 4--6, 
Masson, Paris, 1981. 
\bibitem{Bri1}
\emph{E.~Brieskorn}, Sur les groupes de tresses 
[d'apr\`es V. I. Arnol'd]. (French) S\'eminaire Bourbaki, 24\`eme
   ann\'ee (1971/1972), Exp. No. 401, pp. 21--44. Lecture Notes in Math., 
Vol. 317, Springer, Berlin, 1973. 
\bibitem{Brin}
\emph{M.~Brin},   The Algebra of Strand Splitting. I. A Braided Version of Thompson's 
Group V. Arxiv math.GR/0406042
 \bibitem{BMR}
\emph{M.~Brou\'e, G.~Malle, R.~Rouquier},
Complex reflection groups, braid groups, Hecke algebras.
J. Reine Angew. Math. 500, 127-190 (1998).        
\bibitem{CFP}
\emph{J.~ W.~Cannon, W.~J.~Floyd, W.~R.~Parry}, Introductory notes on Richard 
Thompson's groups.  Enseign. Math. (2)  42  (1996),  no. 3-4, 215--256. 
\bibitem{Ch}
\emph{V.~Chaynikov}, Word and conjugacy problems for the 
singular braid monoids. Comm. Algebra.~4 (2006), no.~6,  1981 - 1995.  
\bibitem{CPVW}
\emph{F. R.~Cohen, J. ~Pakianathan, V.~V.~Vershinin, J.~Wu},
 Basis-conjugating automorphisms of a free group and associated Lie algebras.
 Geometry and Topology Monograph 13 (2008), 147--168.
\bibitem{Cor}
\emph{R.~Corran}, A normal form for a class of monoids including the 
singular braid monoids. J. Algebra. 223
   (2000), no. 1, 256--282. 
\bibitem{CM}
\emph{ H.~S.~M.~Coxeter, W.~O.~J.~Moser},
Generators and relations for discrete groups. 3rd ed.
Ergebnisse der Mathematik und ihrer Grenzgebiete. Band 14.
Berlin-Heidelberg-New York: Springer-Verlag. IX, 161 p. (1972). 
\bibitem{Dehn}
\emph{P.~Dehornoy}, Groupes de Garside.  Ann. Sci. 
\'Ecole Norm. Sup. (4)  35  (2002),  no. 2, 267--306.
\bibitem{Dehn2}
\emph{P.~Dehornoy}, The group of parenthesized braids,Adv. Math. 205 (2006), no. 2, 354--409. 
 \bibitem{DP}
\emph{P.~Dehornoy; L.~Paris}, Gaussian groups and Garside groups, two generalisations 
of Artin groups. Proc. London Math. Soc. (3) 79 (1999), no. 3, 569--604.
\bibitem{Del}
\emph{P.~Deligne},
 Les immeubles des groupes de tresses g\'en\'eralis\'es. (French) 
Invent. Math. 17 (1972),   273--302. 
\bibitem{EL}
\emph{D.~Easdown, T.~G.~Lavers},  The inverse braid monoid.  Adv. Math.  186  
(2004),  no. 2, 438--455.
\bibitem{EM}
\emph{E.~El-Rifai, H.~R.~Morton}, Algorithms for positive braids. 
Quart. J. Math. Oxford Ser. (2) 45 (1994), no.
   180, 479--497. 
\bibitem{E_Th}
\emph{D.~B.~A.~Epstein, J.~W.~Cannon, D.~E.~Holt, S.~V.~F.Levy, 
M.~S.~Paterson, W.~P.~Thurston}, Word processing in groups. 
Jones and Bartlett Publishers, Boston, MA, 1992. xii+330 pp.  
\bibitem{EF}
\emph{B.~Everitt, J.~Fountain},
Partial mirror symmetry I: reflection monoids. 
Adv. Math. 223 (2010), no. 5, 1782-1814,
\bibitem{FaV}
\emph{E.~Fadell, J.~Van Buskirk}, The braid groups of $E\sp{2}$ and $S\sp{2}$. 
Duke Math. J. 29 1962,    243--257.
\bibitem{FarS}
\emph{D.~Farley, L.~ Sabalka}, Discrete Morse theory and graph braid groups.  Algebr. Geom. Topol. 
 5  (2005), 1075--1109.
\bibitem{FKR}
\emph{R.~Fenn, E.~Keyman, C.~Rourke}, The singular braid monoid embeds 
in a group. J. Knot Theory
   Ramifications 7 (1998), no. 7, 881--892.
\bibitem{FRR2}
\emph{R.~Fenn, R.~Rim\'anyi, C.~Rourke}, 
The braid-permutation group. 
Topology 36 (1997), no. 1,
   123--135. 
   \bibitem{g} 
   \emph{S.~Galatius}, Stable homology of automorphism groups of free groups,  math.AT/0610216.
\bibitem{Gar}
\emph{F.~A.~Garside}, 
The braid group and other groups,
Quart. J. Math. Oxford Ser. 1969, 20,
235--254.
\bibitem{Ge1}
\emph{B.~Gemein}, Singular braids and Markov's theorem. J. Knot Theory Ramifications 6 (1997), no. 4, 441--454.
\bibitem{Ghr}
\emph{R.~Ghrist},  Configuration spaces and braid groups on graphs in robotics.  
Knots, braids, and mapping class groups---papers dedicated to 
Joan S. Birman (New York, 1998),  29--40, AMS/IP Stud. Adv. Math., 24, 
Amer. Math. Soc., Providence, RI, 2001.
\bibitem{Gil}
\emph{N.~D.~Gilbert}, Presentations of the inverse braid monoid.  
J. Knot Theory Ramifications  15  (2006),  no. 5, 571--588.
\bibitem{Gol}
\emph{D.~L.~Goldsmith}, 
The theory of motion groups.
Michigan Math. J. 28 (1981), no. 1, 3--17.    
\bibitem{gon3} 
\emph{J.~Gonz\'alez-Meneses},  
Presentations for the monoids of singular braids on closed 
  surfaces, Comm. Algebra,  {30} (2002), 2829-2836.  
\bibitem{GrSer}
\emph{P.~Greenberg, V.~Sergiescu}, An acyclic extension of the braid group. 
 Comment. Math. Helv.  66  (1991),  no. 1, 109--138. 
\bibitem{JMM}
\emph{C.~Jensen, J.~McCammond, J.~Meier}, The integral cohomology of the group of loops.  Geom. Topol.  10  (2006), 759--784.  
\bibitem{Joh1}
\emph{D.~L.~Johnson}, Towards a characterization of smooth braids. 
Math. Proc. Cambridge Philos. Soc. 92 (1982), no. 3, 425--427.  
\bibitem{KT}
\emph{C. Kassel, V.~Turaev}, Braid groups.  Graduate Texts in Mathematics, 247. Springer, New York, 2008. xii+340 pp.
\bibitem{ka} 
\emph{N.~Kawazumi}, Cohomological aspects of Magnus expansions, 
math.GT/0505497.
\bibitem{Kl}
\emph{F.~Klein},
Vorlesungen \"uber h\"ohere Geometrie. 3. Aufl., bearbeitet und herausgegeben 
von {\it W. Blaschke}. VIII${}+{}$405 S. Berlin, J. Springer (Die Grundlehren der 
mathematischen Wissenschaften in Einzeldarstellungen Bd.
22) (1926).
\bibitem{Kou}
\emph{}Kourovka notebook: unsolved problems in group theory. Seventh edition. 
1980.  Akad. Nauk SSSR Sibirsk. Otdel., Inst. Mat., Novosibirsk, 1980. 115 pp. (Russian)
\bibitem{Kr}
\emph{D.~Krammer}, A class of Garside groupoid structures on the pure braid group,
Trans. Amer. Math. Soc. 360 (2008), no. 8, 4029–4061.   
\bibitem{km} 
\emph{ S.~Krsti\'c, J.~McCool}, The non-finite presentability of $IA(F_3)$ and
$GL_2(Z[t,t^{-1}])$, Invent. Math. 129 (1997), 595--606.
\bibitem{La}
\emph{S.~Lambropoulou}, Solid torus links and Hecke algebras of ${\rm B}$-type.  
Proceedings of the Conference on Quantum Topology (Manhattan, KS, 1993),  225--245, 
World Sci. Publishing, River Edge, NJ, 1994. 
\bibitem{Law}
\emph{M.~V.~Lawson}, Inverse semigroups. The theory of partial symmetries. World Scientific Publishing Co., Inc., River Edge, NJ, 1998. xiv+411 pp. 
\bibitem{Li1}
\emph{V. Ya. Lin},  Artinian braids and groups and
spaces connected with them,
 Itogi Nauki i Tekhniki
(Algebra, Topologiya, Geometriya) 1979, 17,
159--227 (Russian). English transl. in J.
Soviet Math. {18} (1982) 736--788.
\bibitem{Li3}
\emph{V. Ya. Lin}, Braids and Permutations,
Arxiv: math.GR/0404528
\bibitem{magnus} W.~Magnus, \emph{\"Uber $n$-dimensionale
Gittertransformationen}, Acta Math., Vol. 64 (1934), 353-367.
\bibitem{mks} \emph{W.~Magnus, A.~Karass, D.~Solitar}, Combinatorial Group Theory, Wiley, 1966.
\bibitem{Mar2}
\emph{A. A. Markoff},  Foundations of the Algebraic Theory of
Tresses, Trudy Mat. Inst. Steklova, No~16, 1945 (Russian, English
summary).
\bibitem{mc} \emph{J.~McCool}, On basis-conjugating automorphisms of free
groups, Canadian J. Math., vol.~38, {12}(1986), 1525-1529.
\bibitem{Mo}
\emph{E.~H.~Moore},  On the reciprocal of the general algebraic matrix. Bull. Amer. Math. Soc. 26,
(1920), 394-395.
\bibitem{v_N}
\emph{J.~von~Neumann}, 
On regular rings. 
Proc. Natl. Acad. Sci. USA 22, 707-713 (1936).
\bibitem{N} 
\emph{J.~Nielsen}, \"Uber die Isomorphismen unendlicher
Gruppen ohne Relation, (German) Math. Ann. 79 (1918), no. 3,
269--272.
\bibitem{orlikterao} 
\emph{P. Orlik, H. Terao},
 Arrangements of hyperplanes,
Grundlehren der mathematischen Wissenschaften {300},
Springer-Verlag, 1992.
\bibitem{Pen}
\emph{R.~Penrose},   A generalized inverse for matrices. Proc. Camb. Phil. Soc. 51,
(1955), 406-413.
\bibitem{Pet}
\emph{M.~Petrich}, Inverse semigroups. Pure and Applied Mathematics (New York). 
A Wiley-Interscience Publication. John Wiley \& Sons, Inc., New York, 1984. x+674 pp.
\bibitem{pe} 
\emph{A.~Pettet}, The Johnson homomorphism and the second cohomology of $IA_n$,  Algebr. Geom. Topol.  5  (2005), 725--740 .
\bibitem{Po}
\emph{L.~M.~Popova},  Defining relations of a semigroup of partial 
endomorphisms of a finite linearly ordered set. (Russian)  Leningrad. 
Gos. Ped. Inst. U\v cen. Zap.  238  1962 78--88. 
\bibitem{s} \emph{T.~Sakasai}, The Johnson homomorphism and the third rational
cohomology group of the Torelli group,  Topology Appl.  148  (2005),  no. 1-3, 83--111.
\bibitem{sa} 
\emph{T.~Satoh}, The abelianization of the congruence
IA-automorphism group of a free group, Math. Proc. Camb. Phil. Soc.
 142  (2007),  no. 2, 239--248.
\bibitem{Sav2}
\emph{A.~G.~Savushkina}, On a group of conjugating automorphisms of a free group. 
(Russian) Mat. Zametki 60 (1996), no. 1, 92--108, 159; translation in Math. 
Notes 60 (1996), no. 1-2, 68--80 (1997).
\bibitem{Sc}
\emph{G.~P.~Scott}, Braid groups and the group of homeomorphisms of a 
surface. Proc. Cambridge Philos. Soc. 68, 1970,  605--617.
\bibitem{Ser}
\emph{V. Sergiescu}, Graphes planaires et
pr\'esentations des groupes de tresses,
Math. Z. 1993, 214,  477--490.
\bibitem{ShT}
\emph{G.~C.~Shephard, J.~A.~Todd},
Finite unitary reflection groups.
Can. J. Math. 6, 274-304 (1954).
\bibitem{T}
\emph{J.~Tits},  Normalisateurs de tores. I. Groupes de Coxeter \'etendus.   
J. Algebra~4  (1966) 96--116.
\bibitem{Va}
\emph{V.~A.~Vassiliev}, Complements of discriminants of smooth maps: 
topology and applications. Translations of Mathematical Monographs, 98. 
American Mathematical Society, Providence, RI, 1992.
   vi+208 pp.
\bibitem{Ve0}
\emph{V.~V.~Vershinin},
Thom spectra of generalized braid groups. Preprint No 95/02-2. Universit\'e 
de Nantes. 1995.
\bibitem{Ve1}
\emph{V.~V.~Vershinin}, On braid groups in handlebodies. 
Sib. Math. J. 39, No.4, 645-654 (1998); translation from 
Sib. Mat. Zh. 39, No.4, 755-764 (1998).
\bibitem{Ve2}
\emph{V.~V.~Vershinin},  On homological properties of singular braids. 
Trans. Amer. Math. Soc. 350 (1998), no. 6,
   2431--2455.
\bibitem{Ve3}
\emph{V.~V.~Vershinin},
Homology of braid groups and their generalizations. Knot theory (Warsaw, 1995), 421–446, Banach Center Publ., 42, Polish Acad. Sci., Warsaw, 1998
\bibitem{Ve4}
\emph{V.~V.~Vershinin}, Generalizations of braids from a homological 
point of view. 
Sib. Adv. Math. 9, No.2, 109-139 (1999). 
\bibitem{Ve6}
\emph{V.~V.~Vershinin}, Braid groups and loop spaces. 
Russ. Math. Surv. 54, No.2, 273-350 (1999); translation from 
Usp. Mat. Nauk 54, No.2, 3-84 (1999).
\bibitem{Ve8}
\emph{V.~V.~Vershinin}, On homology of virtual braids and Burau representation. 
Knots in Hellas '98, Vol. 3 (Delphi). J.~Knot Theory Ramifications 10 (2001), 
no. 5, 795--812.
\bibitem{Ve9}
\emph{V.~V.~Vershinin}, On  presentations of generalizations of braids with few
generators,  Fundam. Prikl. Mat. Vol. 11, No 4, 2005. 23-32.
\bibitem{Ve9_5}
\emph{V.~V.~Vershinin}, Braids, their properties and generalizations. Handbook of algebra. Vol. 4, 427–465, Handb. Algebr., 4, Elsevier/North-Holland, Amsterdam, 2006. 
\bibitem{Ve10} 
\emph{V.~V.~Vershinin}, 
 On the singular braid monoid, (Russian) Algebra i Analiz 21 (2009), no. 5, 19--36; translation in St. Petersburg Math. J. 21 (2010), no. 5, 693-704.
\bibitem{Ve11} 
\emph{V.~V.~Vershinin},  On inverse braid and reflection monoids of type $B$. (Russian) Sibirsk. Mat. Zh. 50 (2009), no. 5, 1010--1015; translation in Sib. Math. J. 50 (2009), no. 5, 798-802. 
 \bibitem{Wag} 
\emph{V.~V.~Wagner}, 
 Generalized groups. (Russian)  Doklady Akad. Nauk SSSR (N.S.)  84,  (1952). 1119--1122.  
\bibitem{Za1}
\emph{O.~Zariski}, On the Poincar\'e
group of rational plane curves,
Am. J. Math.
 1936, 58,  607-619. 
 \end{thebibliography}
\end{document}